# MULTIDIMENSIONAL TRIMMING BASED ON PROJECTION DEPTH[1]

By Yijun Zuo

*Michigan State University*

As estimators of location parameters, univariate trimmed means are well known for their robustness and efficiency. They can serve as robust alternatives to the sample mean while possessing high efficiencies at normal as well as heavy-tailed models. This paper introduces multidimensional trimmed means based on projection depth induced regions. Robustness of these depth trimmed means is investigated in terms of the influence function and finite sample breakdown point. The influence function captures the local robustness whereas the breakdown point measures the global robustness of estimators. It is found that the projection depth trimmed means are highly robust locally as well as globally. Asymptotics of the depth trimmed means are investigated via those of the directional radius of the depth induced regions. The strong consistency, asymptotic representation and limiting distribution of the depth trimmed means are obtained. Relative to the mean and other leading competitors, the depth trimmed means are highly efficient at normal or symmetric models and overwhelmingly more efficient when these models are contaminated. Simulation studies confirm the validity of the asymptotic efficiency results at finite samples.

**1. Introduction.** The sample mean is a very standard estimator of the "center" of a given data set and possesses many desirable properties. Indeed, it is the most efficient estimator at normal models. It, however, is notorious for being extremely sensitive to unusual observations (outliers) and heavy-tailed distributions. Indeed, the mean possesses the lowest breakdown point. To be more robust, the sample median is employed. It has the best breakdown point among all reasonable location estimators. The median, however,

Received November 2004; revised November 2005.
[1]Supported in part by NSF Grants DMS-00-71976 and DMS-02-34078.

is not efficient at normal and other light-tailed distributions. Realizing the drawbacks of the mean and the median and motivated by robustness and efficiency considerations, Tukey in 1948 introduced trimmed means in real data analysis [28]. These estimators can strike a desirable balance between robustness and efficiency and serve as compromises between the two extremes—the mean and the median. Despite numerous competitors introduced since 1948, the robustness and efficiency advantages keep the trimmed mean as the most prevailing estimator of location parameters (see, e.g., [2, 3, 11, 27, 30]).

Data from the real world, however, are often multidimensional and contain "outliers" or "heavy tails." The outliers in high dimensions are far more difficult to detect or identify than in the univariate case since it is often difficult to plot the data and the outliers are not always in the single coordinates. A good sample of a real data set of the latter case is given on page 57 of [23]. A robust procedure such as multidimensional trimming that can automatically detect the outliers or heavy tails is thus desirable. The task of trimming in high dimensions, however, turns out to be nontrivial, for there is no natural order principle in high dimensions. On the other hand, data depth has shown to be a promising tool for providing a center-outward ordering of multidimensional observations; see [14, 29, 37], for example. Points deep inside a data cloud get high depth and those on the outskirts get lower depth. With a depth induced ordering, it becomes quite straightforward to define multivariate trimmed means. Indeed, examples are given in [4, 5, 18, 19, 21, 32], all for Tukey halfspace depth trimming; in [14] and [6] for Liu simplicial depth trimming; and in [34] for general depth trimming (see Section 6 for a detailed discussion).

A natural question raised for depth induced multidimensional trimmed means is: Do they share the same robustness and efficiency advantages of their univariate counterparts over the sample mean? No answer has been given in the literature. Indeed, except for a very few sporadic discussions, very little attention has been paid to the depth based multivariate trimmed means and little is known about their robustness and efficiency. To answer the aforementioned question and to shed light on the robustness and the efficiency aspects of a class of depth trimmed means, the projection depth trimmed means, is the *objective* of this article. Although the paper focuses on *projection depth* trimmed means, the technical approaches are applicable to *other* (such as halfspace) depth trimmed means and covariance matrices as well. Motivation for selecting projection depth is addressed in Section 6.

The paper investigates the local as well as the global robustness of the depth trimmed means via the influence function and breakdown point, respectively. Deriving the influence function of the depth trimmed means is exceptionally involved. The difficulty lies in handling the distribution-dependent depth trimming region.

To investigate the large sample behavior (such as asymptotic relative efficiency) of the sample projection depth trimmed means, we have to establish their limiting distributions. The trimming nature of the estimators makes the study of the asymptotics very challenging. Standard asymptotic theory falls short of the goal. Indeed, even establishing the limiting distribution of the regular univariate trimmed means is not as straightforward as one might imagine. One misconception about this is that the task should be similar to or not much more challenging than the one for the sample mean. In fact, the limiting distribution of the regular trimmed means which were introduced as early as 1948 by Tukey (or perhaps earlier) was not established until 1965 by Bickel. Classical textbooks today still do not prove the limiting distribution and only point out the ad hoc proof of Bickel [2] or Stigler [26] without details. Another misconception about the limiting distribution is that it just follows in a straightforward fashion after one derives the influence function. This actually is not always the case (as shown here in this paper and elsewhere). The challenging task of establishing limiting distributions in this paper for the multidimensional depth trimmed means is accomplished by utilizing empirical process theory (see [32] or [20]).

The paper shows that the projection depth trimmed means (with robust choices of univariate location and scale measures) are highly robust locally (with bounded influence functions) and globally (with the best breakdown point among affine equivariant competitors), as well as highly efficient relative to the mean (and depth medians) at normal and heavy-tailed models. The latter is especially true when the models are slightly contaminated. Findings in the paper indicate that the projection depth trimmed means represent very favorable choices, among the leading competitors, of robust and efficient location estimators for multivariate data.

The rest of the paper is organized as follows. Section 2 introduces projection depth induced regions and trimmed means and discusses some fundamental properties. Section 3 is devoted to the study of the local (the influence function) as well as the global (the finite sample breakdown point) robustness of the depth trimmed estimators. Asymptotic representations and asymptotics are established in Section 4. Section 5 addresses the efficiency issue of the projection depth trimmed means. Concluding remarks in Section 6 end the main body of the paper. Selected proofs of main results and auxiliary lemmas are reserved for the Appendix.

## 2. Projection depth regions and trimmed means.

2.1. *Projection depth functions and regions.* Let $\mu$ and $\sigma$ be univariate location and scale measures of distributions. Typical examples of $\mu$ and $\sigma$ include the pair mean and standard deviation (SD) and the pair median (Med) and median absolute deviations (MAD). Define the *outlyingness* of $x \in \mathbb{R}^d$ with respect to (w.r.t.) the distribution $F$ of $X$ in $\mathbb{R}^d$ ($d \geq 1$) as ([4] and [25])

$$O(x, F) = \sup_{u \in S^{d-1}} |g(x, u, F)|, \tag{1}$$

where $S^{d-1} = \{u : \|u\| = 1\}$, $g(x, u, F) = (u'x - \mu(F_u))/\sigma(F_u)$ is the "generalized standard deviation" of $u'x$ w.r.t. $F_u$ and $F_u$ is the distribution of $u'X$. If $u'x - \mu(F_u) = \sigma(F_u) = 0$, we define the generalized standard deviation $g(x, u, F) = 0$. The *projection depth* of $x \in \mathbb{R}^d$ w.r.t. the given $F$, $PD(x, F)$, is then defined as

$$PD(x, F) = 1/(1 + O(x, F)). \tag{2}$$

Sample versions of $g(x, u, F)$, $O(x, F)$ and $PD(x, F)$ are obtained by replacing $F$ with its empirical version $F_n$. With $\mu$ and $\sigma$ being the Med and the MAD, respectively, Liu [15] first suggested the use of $1/(1 + O(x, F_n))$ as a depth function. Zuo and Serfling [37] defined and studied (2) with $(\mu, \sigma) = $ (Med, MAD). Since $PD$ depends on the choice of $(\mu, \sigma)$, a further study with general $\mu$ and $\sigma$ is carried out in [33]. It turns out $PD$ possesses desirable properties for depth functions (see [37]). For example, it is affine invariant, maximized at the center of a symmetric distribution, monotonically decreasing when a point moves along a ray stemming from the deepest point, and vanishes at infinity. For motivation, examples and other related discussions of (2), see [33].

For any $0 < \alpha < \alpha^* = \sup_{x \in \mathbb{R}^d} PD(x, F) \leq 1$, the $\alpha$*th projection depth region* is

$$PD^\alpha(F) = \{x : PD(x, F) \geq \alpha\}. \tag{3}$$

It is a multivariate analogue of the univariate $\alpha$th quantile region $[F^{-1}(\alpha), F^{-1}(1 - \alpha)]$. The set $\{x : PD(x, F) = \alpha\}$ is called the $\alpha$*th projection depth contour*, which is the boundary $\partial PD^\alpha(F)$ of $PD^\alpha(F)$ under some conditions (see [33]). Structural properties and examples of projection depth regions and contours are discussed in [33]. Note that $\alpha$ in (3) can also be determined by the probability content of the resulting region. For example, define $\alpha(\lambda) = \sup\{\alpha : P_X(x : PD(x, F) \geq \alpha) \geq \lambda\}$; then $P_X(PD^{\alpha(\lambda)}(F)) = \lambda$ for a smooth distribution function $F$. A sample version of $PD^\alpha(F)$, $PD_n^\alpha$, is obtained by replacing $F$ with its empirical version $F_n$.

We assume throughout that $\mu(F_{sY+c}) = s\mu(F_Y) + c$ and $\sigma(F_{sY+c}) = |s|\sigma(F_Y)$ (*affine equivariance*) for any scalars $s$ and $c$ and random variable $Y \in \mathbb{R}^1$, and that

(C0) $\sup_{u \in S^{d-1}} \mu(F_u) < \infty$, $0 < \inf_{u \in S^{d-1}} \sigma(F_u) \leq \sup_{u \in S^{d-1}} \sigma(F_u) < \infty$.

This holds for typical location and scale functionals; see Remark 2.4 of [33]. It follows that $PD^\alpha(F)$ is compact and has a nonempty interior that contains the maximum depth point $\theta$ with $PD(\theta, F) = \alpha^*$ (Theorems 2.2 and 2.3 of [33]). By the affine invariance of the projection depth functions, we can assume without loss of generality that $\theta = 0 \in \mathbb{R}^d$ in our following discussion. The depth region $PD^\alpha(F)$ can then be characterized by the "directional radius functional" $R^\alpha(u, F)$,

(4) $\qquad R^\alpha(u, F) = \sup\{r \geq 0 \colon ru \in PD^\alpha(F)\} \qquad \forall u \in S^{d-1}$,

which is the same as $\inf\{r \geq 0 \colon ru \notin PD^\alpha(F)\}$. For simplicity, we sometimes write $R(u, F)$ or $R(u)$ for $R^\alpha(u, F)$ and $R_n(u)$ for $R^\alpha(u, F_n)$ for fixed $\alpha$ and $F$.

2.2. *Projection depth trimmed means and fundamental properties.* With depth regions, one can define the $\alpha$th *projection depth trimmed mean* (PTM) by

(5) $$PTM^\alpha(F) = \int_{PD^\alpha(F)} w(PD(x, F)) x \, dF(x) \Big/ \int_{PD^\alpha(F)} w(PD(x, F)) \, dF(x),$$

where $w(\cdot)$ is a suitable (bounded) weight function on $[0, 1]$ such that the denominator is nonzero. The latter is true for typical nonzero $w(\cdot)$. Note that the numerator is bounded since $PD^\alpha$ is; see Theorem 2.3 of [33]. Thus $PTM^\alpha(F)$ is well defined. Again we may suppress $\alpha$ and (or) $F$ in $PTM^\alpha(F)$ for convenience.

When $w$ is a (nonzero) constant, (5) gives equal nonzero weight to each point within the depth region $PD^\alpha(F)$, and zero weight to any point outside the region. Thus we have exactly the same (0–1) weighting scheme as that of the regular univariate trimmed mean. Two such $PTM^\alpha(F_n)$'s with $w = c > 0$ are illustrated in Figure 1 with a bivariate standard normal sample of size 900 and $\alpha = 0$ and 0.36. To treat a broader class of multidimensional trimmed means, in our following discussion $w$ is allowed to be any suitable nonconstant function, though.

On the other hand, it is noteworthy that in the degenerate one-dimensional case (with a nonzero constant $w$), (5) yields a *new* type of trimmed mean that is different from the regular one. The difference lies in the trimming scheme. For example, at the sample level the regular trimming is based on the ranks of sample points whereas (5) is based on the values of the generalized standard deviations. The latter can lead to more robust and efficient estimators (see Sections 3.3 and 5).

It can be seen that $PTM^\alpha(\cdot)$ is affine equivariant, that is, $PTM^\alpha(F_{AX+b}) = A(PTM^\alpha(F_X)) + b$ for any nonsingular $d \times d$ matrix $A$ and $b \in \mathbb{R}^d$, since $PD(x, F)$ is affine invariant. Hence $PTM^\alpha$ does not depend on the underlying coordinate system or measurement scale. If $X \sim F$ is centrally symmetric about $\theta \in \mathbb{R}^d$ [i.e., $\pm(X - \theta)$ have the same distribution], then $PTM^\alpha(F)$ is Fisher consistent about $\theta$, that is, $PTM^\alpha(F) = \theta$, and $PTM^\alpha(F_n)$ is also centrally symmetric about $\theta$ since $PTM^\alpha(AX_1 + b, \ldots, AX_n + b) = APTM^\alpha(X_a, \ldots, X_n) + b$ for any nonsingular matrix $A$ and $b \in \mathbb{R}^d$. The latter also implies that $PTM^\alpha(F_n)$ is unbiased for $\theta$.

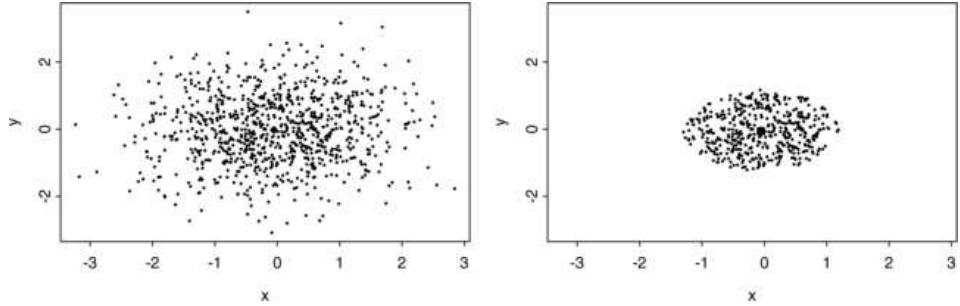

FIG. 1. $PTM^\alpha(F_n)$ based on a $N(0, I_2)$ sample of size 900. Left: $\alpha = 0$ and $PTM^\alpha(F_n) = (-0.05798, -0.02476)$. Right: $\alpha = 0.36$ and $PTM^\alpha(F_n) = (-0.05571, -0.05848)$.

**3. Robustness.** Robustness is a fundamental issue in statistics. It has long been recognized as a principal performance criterion for statistical procedures. We address the *local* and the *global* robustness of depth trimmed means in this section.

One popular (qualitative) robustness measure of a statistical procedure is its influence function. Let $F$ be a given distribution, let $\delta_x$ be the point-mass probability distribution at a fixed point $x \in \mathbb{R}^d$ and let $F(\varepsilon, \delta_x) = (1 - \varepsilon)F + \varepsilon\delta_x$, $\varepsilon \in [0, 1]$, be the point-mass contaminated distribution. The *influence function* (IF) of a statistical functional $T$ at $x \in \mathbb{R}^d$ for the given $F$ is defined as [10]

(6) $$IF(x; T, F) = \lim_{\varepsilon \to 0^+} (T(F(\varepsilon, \delta_x)) - T(F))/\varepsilon,$$

which describes the relative effect (influence) on $T$ of an infinitesimal point-mass contamination at $x$, and captures the local robustness of $T$. A functional with a bounded influence function thus is robust and desirable. The supremum norm of $IF(x; T, F)$ is called the *gross error sensitivity* of $T$ at $F$ [10],

(7) $$GRE(T, F) = \sup_{x \in \mathbb{R}^d} \|IF(x; T, F)\|,$$

the maximum relative effect on $T$ of an infinitesimal point-mass contamination.

It is well known that the mean functional has an unbounded influence function whereas that of the regular univariate trimmed mean functional is bounded; see [24], for example. The natural concern now is whether the influence function of the projection depth trimmed mean functional is bounded.

Note that the integral region in the definition of $PTM^\alpha(F)$ is a functional of $F$. An infinitesimal point-mass contamination hence affects this region. The derivation of the influence function of $PTM^\alpha(F)$ thus becomes challenging. Our strategy to attack the problem is "divide and conquer": to work out the influence function of the projection depth region first and then the influence function of the projection depth region induced trimmed mean functional based on the preliminary results.

3.1. *Influence function of depth region.* Here we establish the influence function of $R^\alpha(u, F)$. Denote by $F_u(\varepsilon, \delta_x)$ the projected distribution of $F(\varepsilon, \delta_x)$ to a unit vector $u$. Then $F_u(\varepsilon, \delta_x) = (1-\varepsilon)F_u + \varepsilon \delta_{u'x}$. For simplicity we sometimes write $F_\varepsilon$ and $F_{\varepsilon u}$ for $F(\varepsilon, \delta_x)$ and $F_u(\varepsilon, \delta_x)$, respectively, for the fixed $x \in \mathbb{R}^d$. We need the following itemized conditions. Denote by $o_x(1)$ a quantity that may depend on a given point $x \in \mathbb{R}^d$ but approaches 0 as $\varepsilon \to 0$ for the fixed $x$.

(C1) $\mu(\cdot)$ and $\sigma(\cdot)$ at $F_u$ and $F_{\varepsilon u}$ are continuous in $u \in S^{d-1}$ and $\sigma(F_u) > 0$,

(C2) $|\mu(F_{\varepsilon u}) - \mu(F_u)| = o_x(1)$, $|\sigma(F_{\varepsilon u}) - \sigma(F_u)| = o_x(1)$ uniformly in $u \in S^{d-1}$,

(C3) $\frac{\mu(F_u(\varepsilon, \delta_x)) - \mu(F_u)}{\varepsilon} = IF(u'x; \mu, F_u) + o_x(1)$, $\frac{\sigma(F_u(\varepsilon, \delta_x)) - \sigma(F_u)}{\varepsilon} = IF(u'x; \sigma, F_u) + o_x(1)$ uniformly in $u \in S^{d-1}$ for fixed $x \in \mathbb{R}^d$.

Conditions (C1)–(C3) hold for smooth $M$-estimators of location and scale (and also for the Med and MAD); see [[12], page 136] and [[33], page 1468]. Note that (C0)–(C3) are connected (nested) in the sense that (C1) implies (C0), (C3) implies (C2) if $IF(u'x; \mu, F_u)$ and $IF(u'x; \sigma, F_u)$ are bounded in $u$, and (C2) holds when (C1) holds and $\mu(F_{\varepsilon u}) - \mu(F_u) = o_x(1)$ and $\sigma(F_{\varepsilon u}) - \sigma(F_u) = o_x(1)$ for any $u = u(\varepsilon) \to u_0$. The latter holds trivially for continuous functionals $\mu(\cdot)$ and $\sigma(\cdot)$, that is, $\mu(G) \to \mu(F)$ and $\sigma(G) \to \sigma(F)$ as $G$ converges weakly to $F$.

When $\mu(F_u)$ and $\sigma(F_u)$ are continuous in $u$ and $\sigma(F_u) > 0$, there is a unit vector $v(x)$ such that $g(x, v(x), F) = O(x, F)$ for $x \in \mathbb{R}^d$. With $v(x)$ we can drop $\sup_{\|u\|=1}$ in the definition of $O(x, F)$, which greatly facilitates technical treatments. Define

(8) $\qquad \mathcal{U}(x) = \{v(x) : g(x, v(x), F) = O(x, F)\}, \qquad x \in \mathbb{R}^d.$

It is usually a singleton (or a finite set) for continuous $F$ and $x \in \partial PD^\alpha(F)$ (see comments after Theorem 2). Indeed, to construct a counterexample is difficult. In the following we consider the case that $\mathcal{U}$ is a singleton for the sake of convenience.

THEOREM 1. *Assume that $IF(v(y)'x; \mu, F_{v(y)})$ and $IF(v(y)'x; \sigma, F_{v(y)})$ are continuous in $v(y)$ for $y \in \partial PD^\alpha(F)$ with $y/\|y\| \in A \subseteq S^{d-1}$ and $\mathcal{U}(y)$ is a singleton for any $y \in \partial PD^\alpha(F)$. Then under (C1)–(C3) with $\beta(\alpha) = (1-\alpha)/\alpha$,*

$$\frac{R^\alpha(u, F(\varepsilon, \delta_x)) - R^\alpha(u, F)}{\varepsilon}$$
$$= \frac{\beta(\alpha) IF(v(y)'x; \sigma, F_{v(y)}) + IF(v(y)'x; \mu, F_{v(y)})}{u'v(y)} + o_x(1),$$

*uniformly in $u \in A$ with $y = R^\alpha(u, F)u$. The influence function of $R^\alpha(u, F)$ is thus given by the first term on the right-hand side.*

The proof of the theorem, technically very demanding and challenging, is given in the Appendix. The influence function of $R^\alpha(u, F)$ at $x$ is determined by those of $\mu$ and $\sigma$ at $v(y)'x$ for the projected distribution $F_{v(y)}$ with $y = R^\alpha(u, F)u$. Since $u'v(y)$ is bounded below from 0 uniformly in $u$ (shown in the proof of the theorem), $IF(x; R^\alpha(u, F), F)$ is bounded as long as those of $\mu$ and $\sigma$ are bounded for $F_{v(y)}$.

The continuity in $v(y)$ of the influence functions of $\mu$ and $\sigma$ at the point $v(y)'x$ for $F_{v(y)}$ with $y \in \partial PD^\alpha$ and $y/\|y\| \in A$ is important. This and the other conditions in the theorem are met with $A = S^{d-1}$ by typical smooth location and scale measures such as the mean and the standard deviation and other $M$-type location and scale measures (see [12]). They are also met by less smooth ones such as the Med and the MAD for suitable (such as elliptically symmetric) distributions. A random vector $X \sim F$ is elliptically symmetric about $\theta$ if $u'(X - \theta) \stackrel{d}{=} \sqrt{u'\Sigma u} Z$ for $u \in S^{d-1}$, some positive definite matrix $\Sigma$ and some random variable $Z \in \mathbb{R}^1$ with $Z \stackrel{d}{=} -Z$, where "$\stackrel{d}{=}$" stands for "equal in distribution." Denote by $F_{\theta,\Sigma}$ such a distribution $F$. Assume, w.l.o.g., that $\theta = 0$ and $\text{MAD}(Z) = m_0$.

COROLLARY 1. *Let $(\mu, \sigma) = (\text{Med}, \text{MAD})$ and $F = F_{\theta,\Sigma}$ with $Z$ having a density $h_z$ that is continuous and positive in small neighborhoods of $0$ and $m_0$. Then*

$$(R^\alpha(u, F(\varepsilon, \delta_x)) - R^\alpha(u, F))/\varepsilon$$
$$= \left(\frac{\beta(\alpha) \operatorname{sign}(|u'\Sigma^{-1}x| - \|\Sigma^{-1/2}u\|m_0)}{4h_z(m_0)} + \frac{\operatorname{sign}(u'\Sigma^{-1}x)}{2h_z(0)}\right) \Big/ \|\Sigma^{-1/2}u\|$$
$$+ o_x(1),$$

*uniformly in $u \in A$, where $A$ is $S^{d-1}$ if $x = 0$ or consists of all $u \in S^{d-1}$ except those $u$'s with $u'\Sigma^{-1}x = 0$ or $u'\Sigma^{-1}x = \pm\|\Sigma^{-1/2}u\|m_0$. Hence the influence function of $R^\alpha(u, F)$, the first term on the right-hand side, is bounded for $u \in A$.*

By Theorem 1 and Corollary 1, $IF(x; R^\alpha(u,F), F)$ is continuous in $v(y)$ with $y = R^\alpha(u,F)u$ for $u \in A$ and depends on $\alpha$ through $\beta(\alpha)$ only. Its existence and behavior for $u \in S^{d-1} - A$ are of little interest for $IF(x; PTM^\alpha(F), F)$, the ultimate goal of all the discussion in this subsection, and not covered by the above results.

The influence function in Corollary 1 is bounded in $x \in \mathbb{R}^d$ for any $u \in A$. This, however, is not true if we select nonrobust $\mu$ and $\sigma$. For example, if $\mu$ and $\sigma$ are the mean and the standard deviation (SD), then for $u \in A = S^{d-1}$

$$IF(x; R^\alpha(u,F), F)$$
$$= \left(\frac{(1-\alpha)((u'\Sigma^{-1}x)^2 - (\|\Sigma^{-1/2}u\|\sigma_z)^2)}{2\alpha\|\Sigma^{-1/2}u\|\sigma_z} + u'\Sigma^{-1}x\right)\bigg/u'\Sigma^{-1}u,$$

with $\sigma_z^2 = \text{var}(Z)$, which is no longer bounded in $x \in \mathbb{R}^d$. To illustrate graphically this influence function and the one in the corollary, we consider $F = N_2(0, I)$ and $\alpha = 0.2$ for simplicity. By orthogonal equivariance, we can just consider $u_0 = (1,0)'$. The influence functions for (Med, MAD) and (mean, SD) become respectively

$$\text{sign}(|x_1| - c)/f(c) + \text{sign}(x_1)/(2f(0)), \quad 2x_1^2 + x_1 - 2 \quad \text{for } x = (x_1, x_2)',$$

with $c = \Phi^{-1}(3/4)$ and $f$ the density of $N(0,1)$, which are plotted in Figure 2.

Figure 2 indicates that $IF(x; R^\alpha(u,F), F)$ with $(\mu, \sigma) = $ (Med, MAD) is a step function jumping at $|x_1| = 0$ and $c$ and is bounded, whereas with $(\mu, \sigma) = $ (mean, SD) it is continuous everywhere but unbounded.

Equipped with preliminary results in this subsection, we are now in position to pursue the influence function of the projection depth region induced means.

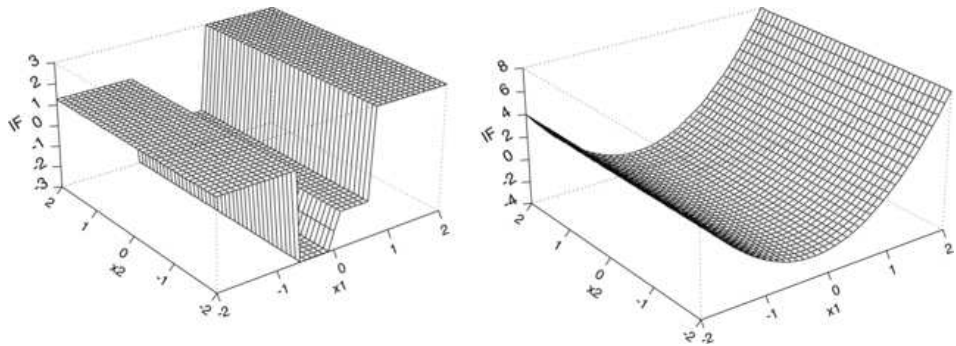

FIG. 2. *The influence functions of $R^\alpha(u,F)$ with $F = N_2(0,I)$, $\alpha = 0.2$ and $u = (1,0)'$. Left: $(\mu, \sigma) = $ (Med, MAD); right: $(\mu, \sigma) = $ (mean, SD).*

3.2. *Influence function of depth trimmed means.* To work out the influence function of the depth trimmed mean functional, we need these conditions:

(C4) $\mathcal{U}(y)$ is a singleton for $y \in B \subseteq PD^\alpha(F)$ with $P_F(PD^\alpha(F) - B) = 0$,

(C5) $IF(u'x; \mu, F_u)$ and $IF(u'x; \sigma, F_u)$ are bounded in $u \in S^{d-1}$ and continuous in $u$ for $u \in \{v(y) : y \in B\}$ with $\int_{S^{d-1} - \{v(y) : y \in B \cap \partial PD^\alpha\}} du = 0$.

In light of the discussion and examples in the last subsection, (C4)–(C5) hold for continuous $F$ and common location and scale functions $\mu$ and $\sigma$ in general. Under these conditions and by virtue of Theorem 1, $IF(x; R^\alpha(u, F), F)$ exists for fixed $x$, $\alpha$ and $F$ and for any $u \in A = \{y/\|y\| : y \in B \cap \partial PD^\alpha\}$ with $\int_{S^{d-1} - A} du = 0$. Now assume that $F^{(1)} = f$ and $w^{(1)}$ exists; then we can define for fixed $\alpha$

$$l_1(x) = \frac{\int_{S^{d-1}} (R(u)u - PTM(F))w(\alpha)f(R(u)u)|J(u, R(u))|IF(x; R(u), F)\,du}{\int_{PD^\alpha(F)} w(PD(y, F))\,dF(y)},$$

$$l_2(x) = \frac{\int_{PD^\alpha(F)} (y - PTM(F))w^{(1)}(PD(y, F))h(x, y)\,dF(y)}{\int_{PD^\alpha(F)} w(PD(y, F))\,dF(y)},$$

$$l_3(x) = \frac{(x - PTM(F))w(PD(x, F))I(x \in PD^\alpha(F))}{\int_{PD^\alpha(F)} w(PD(y, F))\,dF(y)},$$

where

$$h(x, y) = \frac{O(y, F)IF(v(y)'x; \sigma, F_{v(y)}) + IF(v(y)'x; \mu, F_{v(y)})}{\sigma(F_{v(y)})(1 + O(y, F))^2},$$

and $J(u, r)$ is the Jacobian of the transformation from $x \in \mathbb{R}^d$ to $(u, r) \in S^{d-1} \times [0, \infty)$. If we let $x_1 = r\cos\theta_1, \ldots, x_{d-1} = r\sin\theta_1 \sin\theta_2 \cdots \sin\theta_{d-2}\cos\theta_{d-1}$, $x_d = r\sin\theta_1 \cdots \sin\theta_{d-2}\sin\theta_{d-1}$, then $u = x/r$ and $J(u, r) = r^{d-1}\sin^{d-2}\theta_1 \cdots \sin\theta_{d-2}$.

THEOREM 2. *Assume that $F$ has a density $f$ that is continuous in a small neighborhood of $\partial PD^\alpha(F)$ and $w(\cdot)$ is continuously differentiable. Then under (C1)–(C5), $IF(x; PTM^\alpha(F), F) = l_1(x) + l_2(x) + l_3(x)$, which is bounded as long as the influence functions of $\mu$ and $\sigma$ are bounded at $F_u$ for any fixed $u$.*

Note that $\mu$ and $\sigma$ in the theorem can be very general, including mean and SD, Med and MAD, and general $M$-functionals of location and scale. For robust choices, $PTM^\alpha$ has a bounded influence function and hence is (locally) robust.

Note that $\mathcal{U}(y)$ usually is a singleton for $y \in \mathbb{R}^d$ (with the center of symmetry for symmetric $F$ as an exception). For example, if $F = F_{0,\Sigma}$, then $\mathcal{U}(y) = \Sigma^{-1}y/\|\Sigma^{-1}y\|$ for all $y \neq 0$ and *any* affine equivariant $\mu$ and $\sigma$. The condition (C4) in the theorem thus is quite mild. The uniqueness of $v(y)$ ensures a unique limit of the counterpart $v_\varepsilon(y)$ of $v(y)$ as $\varepsilon \to 0$ (see the proof of the theorem). The continuity of $IF(u'x;\mu, F_u)$ and $IF(u'x;\sigma, F_u)$ in $u$ for $u \in \{v(y): y \in B\}$ is sufficient for invoking the result in Theorem 1 and ensures the existence of $h(x,y)$, the influence function of $PD(y,F)$ at point $x$ for any $y \in B$. Conditions in the theorem are usually met by smooth $(\mu, \sigma)$'s such as (mean, SD) and also by less smooth ones such as (Med, MAD).

When $w$ is a nonzero constant (a special yet important case), the influence function $IF(x; PTM^\alpha(F), F)$ becomes $l_1(x) + l_3(x)$ with both terms greatly simplified.

On the other hand, for specific $(\mu, \sigma)$ and $F$ such as (Med, MAD) and $F = F_{\theta,\Sigma}$, the result in the theorem can be concretized. Since it is readily seen that

$$(9) \qquad IF(x; PTM^\alpha, F_{\theta,\Sigma}) = \Sigma^{1/2} IF(\Sigma^{-1/2}(x-\theta); PTM^\alpha, F_{0,I}),$$

we thus will focus on the case $\theta = 0$ and $\Sigma = I$ without loss of generality. We have:

COROLLARY 2. *Let $(\mu, \sigma) = $ (Med, MAD), $F = F_{0,I}$ with density $h_z$ of $Z$ continuous and $> 0$ in small neighborhoods of $0$ and $m_0$, and let $w^{(1)}$ be continuous. Then*

$$IF(x; PTM^\alpha(F), F)$$
$$= \left( \int_{S^{d-1}} c(\alpha)w(\alpha)f(c(\alpha)u)u|J(u,c(\alpha))| IF(x; R(u), F)\, du \right.$$
$$\left. + xc_1/\|x\| + xw(1/(1+\|x\|/m_0))I(\|x\| \leq c(\alpha)) \right) \Big/ c_0,$$

with $c_0 = \int_{\|y\| \leq \beta(\alpha)} w(1/(1+\|y\|))\, dF_0(y)$, $c_1 = \int_{\|y\| \leq \beta(\alpha)} (m_0|y_1|w^{(1)}(1/(1+\|y\|)))/(2h_z(0)(1+\|y\|)^2)\, dF_0(y)$, $y = (y_1,\ldots,y_d)'$, $c(\alpha) = \beta(\alpha)m_0$ and $F_0(y) = F(m_0 y)$.

The most desirable feature of an influence function, the boundedness, is guaranteed by the corollary. This, of course, is no longer true if we select nonrobust $(\mu, \sigma)$ such as (mean, SD). To illustrate this, we consider for simplicity $F = N_2(0, I)$ and a nonzero constant weight function $w$. The influence functions of $PTM^\alpha$ for $(\mu, \sigma) = $ (mean, SD) and (Med, MAD) in this setting at $x = (x_1, x_2)'$ are respectively

$$\left( \int_0^{2\pi} \beta^2(\alpha)g(\beta(\alpha))u(2(u'x)^2 + u'x - 2)\, d\theta \right.$$
$$\left. + xI(\|x\| \leq \beta(\alpha)) \right) \Big/ P_F(\|y\| \leq \beta(\alpha))$$

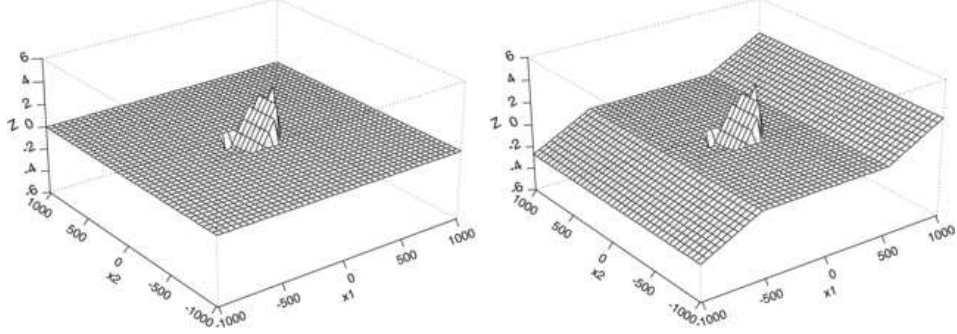

Fig. 3. *The first coordinate of the influence function $IF(x; PTM^\alpha(F), F)$ with $F = N_2(0, I)$ and $\alpha = 0.2$. Left:* $(\mu, \sigma) = (\text{Med}, \text{MAD})$*; right:* $(\mu, \sigma) = (\text{mean}, \text{SD})$.

and

$$\frac{\int_0^{2\pi}(c\beta(\alpha))^2 g(c\beta(\alpha))u(\beta(\alpha)\operatorname{sign}(|u'x|-c)/(4f(c)) + \operatorname{sign}(u'x)/(2f(0)))\,d\theta}{P_F(\|y\| \leq c\beta(\alpha))}$$
$$+ \frac{xI(\|x\| \leq c\beta(\alpha))}{P_F(\|y\| \leq c\beta(\alpha))}$$

with $g(r) = e^{-r^2/2}/(2\pi)$ and $u = (\cos\theta, \sin\theta)'$ (and $c$ defined after Corollary 1), which depend on $\alpha$ through $\beta(\alpha)$ only and are plotted in Figure 3 with $\alpha = 0.2$.

Note that the influence functions in this example are two-dimensional and the figure plots their first coordinates only. The graphs of the second coordinates, however, are the same as the ones in the figure up to an orthogonal transformation.

Both influence functions are continuous except at points $x$ with $\|x\| = c\beta(\alpha)$ or $\beta(\alpha)$. When $\|x\|$ is smaller than these values, the corresponding influence functions behave (roughly) linearly in $x$. The influence of $PTM^\alpha$ with (Med, MAD) is almost zero when $\|x\| > c\beta(\alpha)$. However, in the case with (mean, SD) it becomes unbounded eventually as $\|x\| \to \infty$. All these are reflected clearly in Figure 3.

3.3. *Finite sample breakdown point.* The projection depth trimmed means with robust choices of $\mu$ and $\sigma$ have bounded influence functions and thus are locally robust. This raises the question as to whether they are also globally robust. We now answer this question via the finite sample breakdown point, a notion introduced by Donoho and Huber [7] that has become a prevailing quantitative measure of global robustness of estimators. Roughly speaking, the breakdown point of a location estimator $T$ is the minimum fraction of "bad" (or contaminated) points in a data set that can render $T$ beyond any bound. More precisely, the finite sample breakdown point of $T$ at the sample $X^n = \{X_1, \ldots, X_n\}$ in $\mathbb{R}^d$ is defined as

$$\text{(10)} \qquad \text{BP}(T_n, X^n) = \min\left\{\frac{m}{n} : \sup_{X_m^n} \|T_n(X_m^n) - T_n(X^n)\| = \infty\right\},$$

where $X_m^n$ denotes a contaminated data set resulting from replacing $m$ original points of $X^n$ with $m$ arbitrary points. For a scale estimator $S$, we can calculate its breakdown point by treating $\log S$ as $T$ in the above definition.

Clearly one bad point can ruin the sample mean, hence it has a breakdown point $1/n$, the lowest possible value. The univariate $\alpha$th trimmed mean (trimming $\lfloor \alpha n \rfloor$ data points at both ends of the data) has a breakdown point $(\lfloor \alpha n \rfloor + 1)/n$, which can be much higher than that of the mean. Here $\lfloor \cdot \rfloor$ is the floor function. So the univariate trimmed means can serve as robust alternatives to the sample mean.

For a projection depth trimmed mean, its breakdown point clearly depends on the choice of $(\mu, \sigma)$ in the definition of $PD$. Typical robust choices of $(\mu, \sigma)$ include robust $M$-estimators of location and scale such as (Med, MAD). In the following discussion we first confine attention to the robust choice (Med, MAD) and then comment on the general choices of $(\mu, \sigma)$. We also modify MAD slightly so that the resulting scale measure is less likely to be 0 and consequently the resulting $PD$ trimmed mean has a higher breakdown point. Specifically, we use for $1 \leq k \leq n$

$$\text{MAD}_k = \text{Med}_k\{|x_i - \text{Med}\{x_i\}|\},$$

$$\text{Med}_k\{x_i\} = (x_{(\lfloor (n+k)/2 \rfloor)} + x_{(\lfloor (n+k+1)/2 \rfloor)})/2,$$

where $x_{(1)} \leq \cdots \leq x_{(n)}$ are the ordered values of $x_1, \ldots, x_n$ in $\mathbb{R}^1$. The same idea of modifying MAD to achieve a higher breakdown point for the related estimators has been employed in [31], [9] and [33], for example. Note that when $k = 1$, $\text{MAD}_k$ is just the regular MAD.

For projection depth (or any other depth) trimming, an important issue in practice is how to determine an appropriate value of $\alpha$ so that $PD_n^\alpha \cap X^n$ is not empty and hence $PTM_n^\alpha$ is well defined. It can be shown (based on empirical process theory) that $PD_n^\alpha \cap X^n$ is nonempty almost surely for suitable $\alpha$ under some mild conditions including $P_F(PD^\alpha(F)) > 0$ and *sufficiently* large sample size $n$.

For univariate data, a "pre-data" approach of determining a value of $\alpha$ can be employed in practice. In this case it is not difficult to see that the projection depth of the order statistic $X_{(\lfloor (n+1)/2 \rfloor)}$ is always no less than $1/2$. Hence $PD_n^\alpha \cap X^n$ is nonempty as long as $\alpha \leq 1/2$. For multidimensional data, a "post-data" approach can be adopted. That is, the value of $\alpha$ is data-dependent and determined after $X^n$ becomes available. Since we have to calculate $PD(X_i, F_n)$ anyhow, an appropriate value of $\alpha$ for the trimming can be determined afterward. Or we may select a data-dependent $\alpha$ so that $PD_n^\alpha \cap X^n$ is nonempty, as is done in the following result.

A data set $X^n$ in $\mathbb{R}^d$ $(d \geq 1)$ is *in general position* if there are no more than $d$ sample points contained in any $(d-1)$-dimensional hyperplane. This is true almost surely if the sample is from an absolutely continuous distribution $F$. We have:

THEOREM 3. *Let* $(\mu, \sigma) = (\mathrm{Med}, \mathrm{MAD}_k)$ *with* $k = 1$ *for* $d = 1$ *and* $k = d + 1$ *for* $d > 1$ *and let* $X^n$ *be in general position and* $n \geq d + 1$ *for* $d > 1$. *Then*

$$\mathrm{BP}(PTM^\alpha, X^n) = \begin{cases} \dfrac{\lfloor (n+1)/2 \rfloor}{n}, & d = 1, 0 < \alpha \leq 1/2, \\ \dfrac{\lfloor (n-d+1)/2 \rfloor}{n}, & d > 1, 0 < \alpha \leq \alpha_d, \end{cases}$$

*where the weight* $w(r) > 0$ *for* $0 < r \leq 1$ *and* $\alpha_d := \alpha_d(X^n)$ $(d > 1)$ *satisfying*

$$\frac{1 - \alpha_d(X^n)}{\alpha_d(X^n)} = \max_{\|u\|=1} \frac{\max_{i_1, \ldots, i_{\lfloor (n+2)/2 \rfloor}} \max_{1 \leq k, l \leq \lfloor (n+2)/2 \rfloor} |u'(X_{i_k} - X_{i_l})|}{\min_{i_1, \ldots, i_{d+1}} \max_{1 \leq k, l \leq (d+1)} |u'(X_{i_k} - X_{i_l})|/2},$$
(11)

*and* $i_1, \ldots, i_r$ *are* $r$ *arbitrary distinct integers from the set* $\{1, 2, \ldots, n\}$.

Note that the denominator on the right-hand side of (11) is bounded below from 0 uniformly in $u$ since $X^n$ is in general position. Hence $\alpha_d$ is well defined. It is also seen that $\alpha_d(X^n)$ is *affine invariant*, that is, $\alpha_d(AX^n + b) = \alpha_d(X^n)$ for any nonsingular matrix $A$ and $b \in \mathbb{R}^d$, where $AX^n = \{AX_1, \ldots, AX_n\}$. Thus $PTM^\alpha$ is affine equivariant. Clearly $\alpha_d(X^n) < \min\{\frac{1}{3}, \sup_i PD(X_i, F_n)\}$. Indeed, it is seen that $O(X_i, X_m^n)$ is no greater than the right-hand side of (11) for any original $X_i \in X^n$ and any $0 \leq m \leq \lfloor (n-d+1)/2 \rfloor - 1$. Hence $PTM^\alpha$ for $d > 1$ is well defined.

The main idea of the proof can be briefly explained as follows. The estimator breaks down only if $PD^\alpha(X_m^n)$ is empty or contains points of $X_m^n$ with arbitrarily large norms. This cannot happen unless $\mu$ and (or) $\sigma$ break(s) down. To break down Med, $\lfloor (n+1)/2 \rfloor/n$ contaminating points are needed. With the contaminating points $m = \lfloor (n-d+1)/2 \rfloor$ in $\mathbb{R}^d$ $(d > 1)$, we can force $\mathrm{MAD}_{d+1}(u'X_m^n)$ for the special $u$ to be zero or unbounded. All these can lead to the breakdown of $PTM^\alpha$.

The breakdown point results in the theorem are striking. In $R^1$, $PTM^\alpha$ for any $\alpha \in (0, 1/2]$ achieves the best breakdown point of any translation equivariant location estimators (see [17]). Note that the breakdown point of the regular $\alpha$th trimmed mean is only $(\lfloor \alpha n \rfloor + 1)/n$, which is lower unless $\alpha \approx 0.50$ (which corresponds to the median). The difference in breakdown point between the two types of trimmed means is due to the difference in trimming. In the projection depth trimming case, trimming is done based on the values of the $|X_i - \text{Med}(X^n)|/\text{MAD}(X^n)$'s, while in the regular trimming (equivalent to Tukey halfspace depth trimming) case, it is done based on the ranks of the $X_i$'s.

In $R^d$ $(d > 1)$, the breakdown point, $\lfloor (n-d+1)/2 \rfloor /n$, is also the best (highest) among existing affine equivariant location estimators, with very few exceptions.

Note that the theorem allows one to select very small $\alpha$ values (e.g., 0.05 or 0.10), which then can lead to very high efficiency for $PTM^\alpha$ at (contaminated) normal models (see Section 5), while enjoying very high breakdown point robustness.

For simplicity, the theorem reports only the best breakdown points with the corresponding $d$ and $k$. For general $k$ and $d$, it can be shown that $\text{BP}(PTM^\alpha, X^n)$ is $\lfloor (n-k+2)/2 \rfloor /n$ for $d=1$ and $\min\{\lfloor (n+k+1-2d)/2 \rfloor, \lfloor (n-k+2)/2 \rfloor\}/n$ for $d > 1$. The theorem can be extended for arbitrary $X^n$. In this case, the BP results still hold if $d$ is replaced by $c(X^n)$, the maximum number of sample points contained in any $(d-1)$-dimensional hyperplane. The theorem considers robust choices of $\mu$ and $\sigma$. It can be extended for more general cases. For example, if (mean, SD) is used, then $\text{BP}(PTM^\alpha, X^n)$ is $1/n$. For general $(\mu, \sigma)$, the BP of $PTM^\alpha$ is no less than the minimum of the BPs of $\mu$ and $\sigma$ at $u'X^n$ for arbitrary $u \in S^{d-1}$.

**4. Asymptotics.** This section investigates the large sample behavior of the sample projection depth trimmed means. We focus on the strong consistency and the asymptotic normality of the estimators. To this end, we (have to) characterize the asymptotic behavior of the random convex and compact set $PD^\alpha(F_n)$, the sample projection depth region, via that of the random directional radius $R^\alpha(u, F_n)$.

4.1. *Strong consistency and asymptotic representation of the directional radius.* Denote by $F_{nu}$ the empirical distribution function of $u'X_i$, $i = 1, \ldots, n$, for $u \in S^{d-1}$. The following counterparts of (C1) and (C2) are needed in the sequel:

(C1$'$) $\mu(\cdot)$ and $\sigma(\cdot)$ at $F_u$ and $F_{nu}$ are continuous in $u \in S^{d-1}$ and $\sigma(F_u) > 0$,

(C2$'$) $\sup_{\|u\|=1} |\mu(F_{nu}) - \mu(F_u)| = o(1)$, $\sup_{\|u\|=1} |\sigma(F_{nu}) - \sigma(F_u)| = o(1)$, a.s.,

which hold for common choices of $(\mu, \sigma)$ and a wide range of distributions $F$; see Remark 2.4 of [33] for a detailed account (also see [35]). Indeed, for general $M$-estimators $\mu$ and $\sigma$ including (Med, MAD), (C1$'$) holds for suitable $F$ (see Lemma 5.1 of [33]), which in turn implies (C2$'$) if $\mu(F_{nu_n})$ and $\sigma(F_{nu_n})$ are strongly consistent for $\mu(F_u)$ and $\sigma(F_u)$, respectively, for

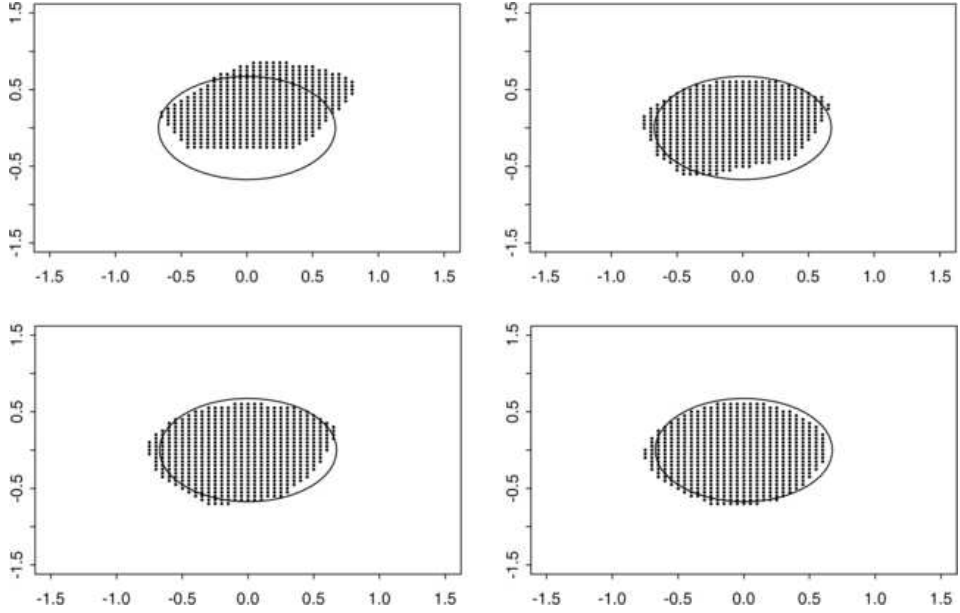

FIG. 4. $R^{0.5}(u, F)$ (*solid circle*) *and* $R^{0.5}(u, F_n)$ (*boundary of the shaded region*) *for* $F = N_2(0, I)$. Upper: *left—*$n = 100$, *right—*$n = 200$. Lower: *left—*$n = 300$, *right—*$n = 900$.

any $u_n \to u \in S^{d-1}$. The latter is true for typical $\mu$ and $\sigma$ since $\mu(G)$ and $\sigma(G)$ are typically continuous in $G$ in the sense that $\mu(G^*)$ and $\sigma(G^*) \to \mu(G)$ and $\sigma(G)$, respectively, whenever $G^*$ becomes close enough to $G$ in distribution (or in Smirnov–Kolmogorov distance) sense (see Example II.1 of [20] for the median functional). We have

THEOREM 4. *Under* (C1′)–(C2′), $\sup_{u \in S^{d-1}} |R^\alpha(u, F_n) - R^\alpha(u, F)| = o(1)$, *a.s.*

The main idea of the proof is as follows. Condition (C1′) insures that for a fixed $x \in \mathbb{R}^d$ there are unit vectors $v(x)$ and $v_n(x)$ such that $O(x, F) = g(x, v(x), F)$ and $O(x, F_n) = g(x, v_n(x), F_n)$ [see (1)]. This result enables us to bound $R^\alpha(u, F_n) - R^\alpha(u, F)$ from above and below for any fixed $u \in S^{d-1}$. Both the upper and the lower bounds are then shown to be $o(1)$ almost surely and uniformly in $u \in S^{d-1}$. A crucial step for this is to show that $x'v(x)$ and $x'v_n(x)$ are bounded below from 0 uniformly for any $x$ on the boundary of $PD^\alpha(F)$ and $PD^\alpha(F_n)$, respectively.

The uniform strong consistency property of $R^\alpha(u, F_n)$ is illustrated in Figure 4. Here $R^\alpha(u, F)$ and $R^\alpha(u, F_n)$ are plotted for $\alpha = 0.5$ and different $n$'s. For simplicity, $F = N_2(0, I)$ is selected. $R^\alpha(u, F)$ then is the circle with radius $\Phi^{-1}(3/4)$. The boundary of $PD^\alpha(F_n)$ is $R^\alpha(u, F_n)$. The uniform strong consistency is clearly demonstrated as $\sup_{\|u\|=1}|R^\alpha(u,F_n) - R^\alpha(u,F)|$ gets smaller when $n$ gets larger.

REMARK 4.1. Under some (stronger) conditions on $F$, $PD^\alpha(F)$ are continuous in Hausdorff distance sense, that is, $\rho(PD^\alpha, PD^{\alpha_0}) \to 0$ as $\alpha \to \alpha_0$, where $\rho(A,B) = \inf\{\varepsilon|\varepsilon > 0, A \subset B^\varepsilon, B \subset A^\varepsilon\}$ and $C^\varepsilon = \{x|\inf\{\|x-y\|:y \in C\} < \varepsilon\}$ (see Theorem 2.5 of [33]). With this continuity of the depth regions, the result in the theorem can be established in a straightforward fashion. For the halfspace depth regions and assuming this continuity, Nolan [19] first obtained the strong consistency result for the radius of the halfspace depth region.

To establish the normality of $R(u, F_n)$, the counterpart of (C3) is needed:

(C3′) The asymptotic representations hold uniformly in $u$:

$$\mu(F_{nu}) - \mu(F_u) = \frac{1}{n}\sum_{i=1}^n f_1(X_i, u) + o_p(n^{-1/2}),$$

$$\sigma(F_{nu}) - \sigma(F_u) = \frac{1}{n}\sum_{i=1}^n f_2(X_i, u) + o_p(n^{-1/2}).$$

The graphs of functions in $\{f_j(\cdot, u): u \in S^{d-1}\}$ form a polynomial discrimination class, $Ef_j(X, u) = 0$ for $u \in S^{d-1}$, $E(\sup_{\|u\|=1} f_j^2(X, u)) < \infty$, for $j = 1$ or 2, and

$$E\left[\sup_{|u_1 - u_2| \leq \delta} |f_j(X, u_1) - f_j(X, u_2)|^2\right] \to 0 \quad \text{as } \delta \to 0, j = 1, 2.$$

For the definition of a class of sets with polynomial discrimination, see [20]. Condition (C3′) holds for general $M$-estimators of $(\mu, \sigma)$ including (Med, MAD) and a wide range of $F$; see [35] for detailed accounts. For example, when $(\mu, \sigma) = (\text{Mean}, \text{SD})$ and $E\|X\|^4$ exists, then $f_1(X, u) = u'(X - EX)$ and $f_2(X, u) = u'((X - EX)(X - EX)' - \text{cov}(X))u/(2\sqrt{u'\text{cov}(X)u})$ and (C3′) holds.

THEOREM 5. *Let $\mathcal{U}(x)$ be a singleton for $x \in \partial PD^\alpha(F)$. Under (C1′)–(C3′),*

$$R^\alpha(u, F_n) - R^\alpha(u, F) = \frac{1}{n}\sum_{i=1}^n k(X_i, R^\alpha(u, F)u) + o_p(n^{-1/2}) \quad \text{a.s.}$$

*uniformly in $u \in S^{d-1}$, where $k(x, y) = (\beta(\alpha)f_2(x, v(y)) + f_1(x, v(y)))/(y_0'v(y))$, for any $y_0 = y/\|y\|$ with $y \neq 0$. Hence*

$$\{\sqrt{n}(R^\alpha(u, F_n) - R^\alpha(u, F)) : u \in S^{d-1}\} \xrightarrow{d} \{Z^\alpha(u) : u \in S^{d-1}\},$$

with $Z^\alpha(u)$ being a zero-mean Gaussian process on the unit sphere with covariance structure $E[k(X, R^\alpha(u_1, F)u_1)k(X, R^\alpha(u_2, F)u_2)]$ for unit vectors $u_1$ and $u_2$.

By virtue of the lower and upper bounds for $R^\alpha(u, F_n) - R^\alpha(u, F)$ established in the proof of Theorem 4 and thanks to empirical process theory (see [32] or [20]), the asymptotic representation for $R^\alpha(u, F_n)$ is obtained after we show that $v_n(R^\alpha(u, F_n)u)$ converges to $v(R^\alpha(u, F)u)$ uniformly in $u \in S^{d-1}$. The directional radius $R^\alpha(u, F_n)$ thus is asymptotically normal for fixed $u \in S^{d-1}$ and also converges as a process to a Gaussian process indexed by $u \in S^{d-1}$. Conditions in the theorem are met by typical $M$-estimators of location and scale and a wide range of distribution functions $F$. For specific $(\mu, \sigma)$, we have specific $k(x, y)$. For example, let $(\mu, \sigma) = (\text{Med}, \text{MAD})$ and $F = F_{0,\Sigma}$; then the following holds.

COROLLARY 3. *Let $(\mu, \sigma) = (\text{Med}, \text{MAD})$, $F = F_{0,\Sigma}$ with the density $h_z$ of $Z$ continuous and $> 0$ in small neighborhoods of $0$ and $m_0$. Then Theorem 5 holds with*

$$f_1(x, u) = \frac{\sqrt{u'\Sigma u}}{h_z(0)}\left(\frac{1}{2} - I(u'x \le 0)\right),$$

$$f_2(x, u) = \frac{\sqrt{u'\Sigma u}}{2h_z(m_0)}\left(\frac{1}{2} - I(|u'x| \le \sqrt{u'\Sigma u}m_0)\right),$$

*and $v(y) = (\Sigma^{-1}y)/\|\Sigma^{-1}y\|$ for any $y \neq 0$.*

The proof of this result is skipped. For related discussion, see Lemma 5.1 of [33] and Lemma 3.2 of [35]. Equipped with the results on $R^\alpha(u, F_n)$, we now are in position to discuss the asymptotics of the depth trimmed means.

4.2. *Strong consistency and asymptotic representation of depth trimmed means.* Strong consistency holds for $PTM^\alpha(F_n)$ under very mild conditions for $\alpha < \alpha^*$.

THEOREM 6. *Let $\mu(F_u)$ and $\sigma(F_u)$ be continuous in $u$ and $\sigma(F_u) > 0$ for $u \in S^{d-1}$ and let $w^{(1)}(\cdot)$ be continuous. Then under (C2′), $PTM(F_n) - PTM(F) = o(1)$, a.s.*

Again the theorem focuses on strong consistency. Other types of consistency can be established accordingly under appropriate versions of (C2′). Note that the weight function in the theorem can be a nonzero constant.

With a standard means, the proof of the theorem seems challenging. The difficulty lies in handling the integral region $PD^\alpha(F_n)$ [or integrand containing $I(x \in PD^\alpha(F_n))$] in (5). The problem becomes less difficult with the help of empirical process theory. The main tool employed in our proof is the generalized Glivenko–Cantelli theorem for a class of measurable functions whose graphs form a class with polynomial discrimination (see II.5 of [20]), or a Glivenko–Cantelli class of measurable functions (see [32]).

We now establish the limiting distribution of $PTM^\alpha(F_n)$ via an asymptotic representation. Assume that $F^{(1)} = f$ and $w^{(1)}$ exist. Replace $IF(x, R^\alpha(u,F), F)$, $IF(v(y)'x, \mu, F_{v(y)})$ and $IF(v(y)'x, \sigma, F_{v(y)})$ with $k(x, R^\alpha(u,F)u)$, $f_1(x, v(y))$ and $f_2(x, v(y))$, respectively, in $l_i(x)$ and $h(x,y)$ and call the resulting functions $\tilde{l}_i(x)$, $i = 1, 2, 3$, and $\tilde{h}(x,y)$, respectively. We have:

THEOREM 7. *Assume that $f$ is continuous in a small neighborhood of $\partial PD^\alpha(F)$, $P_F(\partial PD^\alpha(F)) = 0$ and $w^{(1)}$ is continuous. Then under (C1')–(C3') and (C4)*

$$PTM^\alpha(F_n) - PTM^\alpha(F) = \frac{1}{n}\sum_{i=1}^n (\tilde{l}_1(X_i) + \tilde{l}_2(X_i) + \tilde{l}_3(X_i)) + o_p(1/\sqrt{n}).$$

*Thus $\sqrt{n}(PTM^\alpha(F_n) - PTM^\alpha(F)) \xrightarrow{d} N_d(\mathbf{0}, V)$, where $V = \mathrm{cov}(\tilde{l}_1(X) + \tilde{l}_2(X) + \tilde{l}_3(X))$.*

With standard tools, it seems extremely challenging to establish the asymptotic representation and the normality of $PTM(F_n)$. Thanks to empirical process theory for a Donsker class of functions, especially the asymptotic tightness of the sequence of empirical processes and the asymptotic equicontinuity result and the central limit theorem for the empirical process indexed by a class of functions (see [20] or [32]), we are able to tackle the problem. One key step in our proof is to characterize the complicated integral region $PD^\alpha(G)$ via the directional radius function $R^\alpha(u, G)$ for $G = F$ and $F_n$.

With a nonzero constant $w$, $PTM^\alpha(F_n)$ becomes a depth trimmed mean with equal weight assigned to each sample point within $PD^\alpha(F_n)$. The representation is simplified since $l_2(x)$ vanishes and $l_1(x)$ and $l_3(x)$ also become less complicated.

Conditions in the theorem are met by typical $M$-estimators of location and scale and a wide range of distributions $F$. For example, when $(\mu, \sigma) = (\mathrm{Med}, \mathrm{MAD})$ and $F$ is elliptically symmetric about the origin (assume $\Sigma = I$, w.l.o.g.), we have:

COROLLARY 4. *Let $(\mu, \sigma) = (\mathrm{Med}, \mathrm{MAD})$ and $F = F_{0,I}$ with $F' = f$ being continuous in a small neighborhood of $\|y\| = \beta(\alpha)m_0 = c(\alpha)$ and the density $h_z$ of $Z$ being continuous and positive in small neighborhoods of $0$ and $m_0$. Let $w^{(1)}$ be continuous. Then conditions in Theorem 7 hold and*

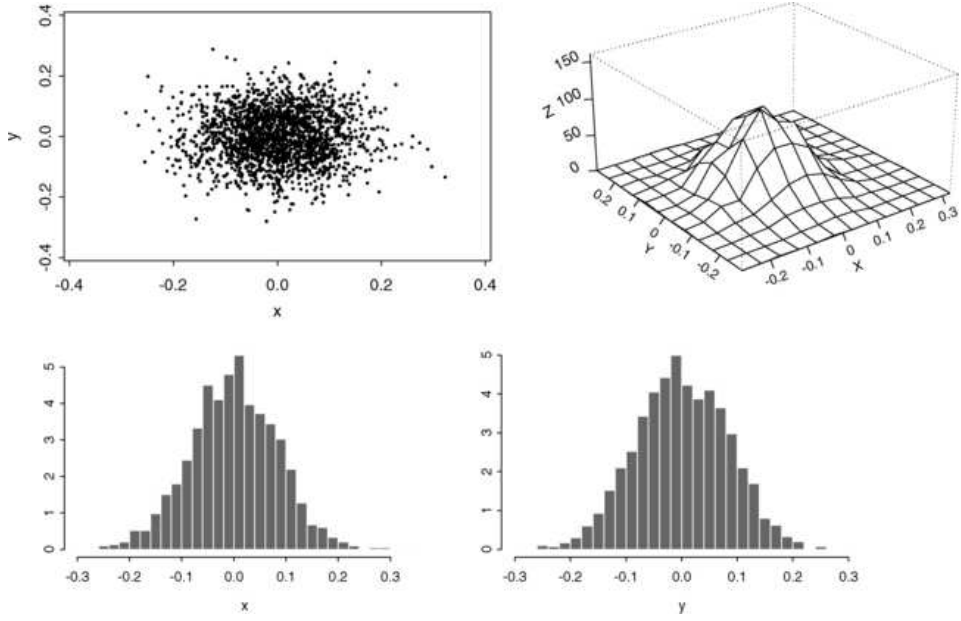

FIG. 5. Upper: *left—plot of* 2000 $PTM_n^\alpha$'s, *right—histogram of* $PTM_n^\alpha$'s. Lower: *coordinate-wise histograms of* $PTM_n^\alpha$'s $[\alpha = 0.36,\ n = 300,\ w = c > 0\ and\ F = N_2(0, I)]$.

$\tilde{l}_1(x) + \tilde{l}_2(x) + \tilde{l}_3(x)$ *is*

$$\left( \int_{S^{d-1}} c(\alpha) w(\alpha) f(c(\alpha)u) |J(u, c(\alpha))| k(x; c(\alpha)u) u\, du \right.$$

$$\left. + \frac{c_1 x}{\|x\|} + xw\left(\frac{1}{1 + \|x\|/m_0}\right) \right) \Big/ c_0,$$

*where* $c_0$ *and* $c_1$ *are defined in Corollary* 2 *and* $k(x; c(\alpha)u) = \beta(\alpha)(\frac{1}{2} - I(|u'x| \le m_0))/(2h_z(m_0)) + (\frac{1}{2} - I(u'x \le 0))/h_z(0)$.

The asymptotic normality in Theorem 7 and Corollary 4 is illustrated in Figure 5. Here 2000 $PTM^\alpha(F_n)$'s are obtained based on $N_2(0, I)$ with $n = 300$ and $\alpha = 0.36$. The two-dimensional histogram indicates a (roughly) normal shape, and so do the one-dimensional histograms of the $x$- and $y$-coordinates of the $PTM^\alpha(F_n)$'s.

**5. Efficiency.** Besides robustness, efficiency is another fundamental issue in statistics. It is also a key performance criterion for any statistical procedure. Section 3 reveals that $PTM^\alpha$ is robust locally and globally for suitable choices of $\mu$ and $\sigma$. A natural question is: Is $PTM^\alpha$ also highly efficient at normal and other models? This section answers the question at both large and finite samples.

5.1. *Large sample relative efficiency.* Consider for simplicity the case that $(\mu,\sigma) = (\text{Med}, \text{MAD})$ and $w = c > 0$. Following the convention in the location setting, assume that $F = F_{\theta,\Sigma}$. By affine equivariance, assume, w.l.o.g., that $\theta = 0$ and $\Sigma = I$. Furthermore, assume that $F' = f$ and the density $h_z$ of $Z$ is continuous and positive at 0 and $m_0$. By Theorems 7 and 5 and Corollaries 4 and 3, we have

COROLLARY 5. *Let $(\mu,\sigma) = (\text{Med}, \text{MAD})$ and $F = F_{0,I}$ meet the conditions in Corollary 4. Let $w = c \neq 0$. Then results in Theorem 4 hold with $\tilde{l}_2(x) = 0$ and*

$$\tilde{l}_1(x) + \tilde{l}_3(x)$$
$$= \frac{\int_{S^{d-1}} c(\alpha) f(c(\alpha)u) |J(u, c(\alpha)| k(x, c(\alpha)u) u \, du + xI(\|x\| \leq c(\alpha))}{P(\|X\| \leq c(\alpha))},$$

*where $k(x, c(\alpha)u) = \beta(\alpha) \operatorname{sign}(|u'x| - m_0)/(4h_z(m_0)) + \operatorname{sign}(u'x)/(2h_z(0))$ and $\sqrt{n}(PTM^\alpha(F_n) - PTM^\alpha(F)) \xrightarrow{d} N_d(\mathbf{0}, V)$, where for $X = (X_1, \ldots, X_d)'$*

$$V = \frac{E(X_1 I(\|X\| \leq c(\alpha)) + \int_{S^{d-1}} c(\alpha) f(c(\alpha)u) |J(u, c(\alpha)| k(X, c(\alpha)u) u_1 \, du)^2}{P^2(\|X\| \leq c(\alpha))}$$
$$\times I_{d \times d}.$$

The key ingredient of the proof of the corollary is repeatedly taking advantage of the symmetry of $F$ in a nontrivial (and clever) manner. The proof, albeit not very challenging technically, is quite involved and hence is skipped.

The explicit form of $V$ greatly facilitates the calculation of the asymptotic relative efficiency of $PTM^\alpha(F_n)$. Note that $EX_1 I(\|X\| \leq c(\alpha)) \operatorname{sign}(|u'X| - m_0) = 0$, which further simplifies the calculation. Call the denominator of $V$ $a$ and the numerator $b$; then $V = b/a I_{d \times d}$. Hence the asymptotic efficiency of $PTM^\alpha(F_n)$ relative to the sample mean is $(a\sigma_z^2)/b$ with $\sigma_z^2 = \operatorname{var}(Z)$. For $X \sim N_d(0, I)$, we have $\sigma_z^2 = 1$, $a = (P(T \leq c^2(\alpha)))^2$ with $T \sim \chi^2(d)$, $f(c(\alpha)u) = g(c(\alpha)) = e^{-c^2(\alpha)/2}/(2\pi)^{d/2}$ and $m_0 = \Phi^{-1}(3/4)$. When $d = 2$, $u = (\cos(\theta), \sin(\theta))$, $a = (1 - e^{-c^2(\alpha)/2})^2$ and

$$b = E\left(X_1 I(\|X\| \leq c(\alpha)) + \int_0^{2\pi} c^2(\alpha) g(c(\alpha)) k(X, c(\alpha)u) \cos(\theta) \, d\theta\right)^2.$$

In Table 1 we list the asymptotic relative efficiency (ARE) results of $PTM^\alpha$ for different $\alpha$'s of the Stahel–Donoho estimator (see [35]) and of the halfspace median (HM) and the projection median (PM) (see [33]) at $N_2(0, I)$.

It is seen that $PTM^\alpha$ is highly efficient for small $\alpha$'s and is much more efficient than some leading competitors. Replacing Med in $PTM^\alpha$ (and PM) with a more efficient one at normal (and other) models, one can improve the efficiencies of $PTM^\alpha$ (and PM); see [33] for discussion related to PM. Our calculations indicate that when the "tail" of $F$ get heavier, $PTM^\alpha$ can gets more efficient than the mean. Furthermore, when $d$ increases, the ARE of $PTM^\alpha$, as expected, increases.

TABLE 1
*ARE of depth trimmed means and medians relative to the mean*

| $PTM^{0.05}$ | $PTM^{0.10}$ | $PTM^{0.15}$ | $PTM^{0.20}$ | SD | PM | HM | Mean |
|---|---|---|---|---|---|---|---|
| 0.9990 | 0.9981 | 0.9927 | 0.8856 | 0.935 | 0.77 | 0.76 | 1.00 |

5.2. *Finite sample relative efficiency.* The comparisons of the relative efficiency results in the last subsection have all been asymptotic, and this raises the question as to whether they are relevant at finite sample practice. Asymptotic results indeed are quite vulnerable to criticism about their practical merits. We now address this issue in this subsection through finite sample Monte Carlo studies.

To see how $PTM^\alpha$ performs in a neighborhood of a normal model, we generate $m = 1000$ samples for different sizes $n$ from the model $(1-\varepsilon)N_2((0,0)', \mathbf{I}) + \varepsilon N_2((\mu,\mu)', \sigma^2 \mathbf{I})$ with $\mu = 10$ and $\sigma = 5$ and $\varepsilon = 0.0$, 0.1 and 0.2. For simplicity, we just consider the case $\alpha = 0.1$. Included in our study are Stahel–Donoho (SD) [35], PM [33] and HM estimators. We assume that all the estimators aim at estimating the known location parameter $\theta = (0,0)' \in \mathbb{R}^2$.

For an estimator $T$ we calculate its "empirical mean squared error" (EMSE) $\sum_{i=1}^m \|T_i - \theta\|^2/m$, where $T_i$ is the estimate based on the $i$th sample. The relative efficiency (RE) of $T$ w.r.t. the mean is obtained by dividing the EMSE of the mean by that of $T$. Here $(\mu, \sigma) = (\text{Med, MAD})$, $w = c \neq 0$. Tables 2–4 list some efficiency results relative to the mean. The entries in parentheses are EMSE$\times 10^3$.

Table 2 reveals that for a perfect $N_2(0, I)$ model $PTM^{0.1}$ is extremely (and the most) efficient. The consistency of RE's with the ARE's confirms the validity of the results in Table 1. The SD estimator is the second most (about 93%) efficient and the PM and HM with roughly the same efficiency are the least efficient ones.

In practice, data more often than not follow a model that is not perfectly normal. Typical examples include contaminated normal (or mixture normal) models. This raises the question of the practical relevance (or robustness) of the results in Table 2. Tables 3 and 4 indicate that $PTM^{0.1}$ has very (most) *robust* EMSE's. Indeed under $\varepsilon = 0.1$ and 0.2, the EMSE's of $PTM^{0.1}$ are still very close to those with $\varepsilon = 0.0$. This robustness increasingly degenerates for SD, PM and HM. The mean has the least robust

TABLE 2
*Finite sample efficiency of $PTM^\alpha$ relative to the sample mean*

| | | | $N_2((0,0)', I)$ | | |
|---|---|---|---|---|---|
| $n$ | $PTM^{0.1}$ | SD | PM | HM | Mean |
| 20 | 0.9843 | 0.9381 | 0.7999 | 0.8053 | 1.0000 |
| | (104.24) | (109.38) | (128.27) | (127.42) | (102.61) |
| 40 | 0.9985 | 0.9298 | 0.7822 | 0.7732 | 1.0000 |
| | (50.560) | (54.299) | (64.546) | (65.296) | (50.485) |
| 60 | 0.9984 | 0.9347 | 0.7675 | 0.7671 | 1.0000 |
| | (32.941) | (35.187) | (42.850) | (42.873) | (32.889) |
| 80 | 1.0000 | 0.9387 | 0.7782 | 0.7762 | 1.0000 |
| | (25.146) | (26.787) | (32.314) | (32.398) | (25.146) |
| 100 | 0.9995 | 0.9338 | 0.7762 | 0.7645 | 1.0000 |
| | (20.014) | (21.421) | (25.770) | (26.061) | (20.003) |

EMSE's. Indeed, the EMSE's of the mean change drastically (enlarged 100 times or more) under the contaminations. With slight departures from normality, all the depth estimators become overwhelmingly more efficient than the mean while $PTM^{0.1}$ performs substantially better than its competitors. The results here for SD, PM and HM are very consistent with those in [35] and [33].

Our simulation studies indicate that the above findings also hold true for other nonnormal (such as $t$, double-exponential and logistic) models. Furthermore, the relative efficiency of $PTM^\alpha$ increases as the dimension $d$ increases.

TABLE 3
*Finite sample efficiency of $PTM^\alpha$ relative to the sample mean*

| | | | $0.90N_2((0,0)', I) + 0.10N_2((10,10)', 25I)$ | | |
|---|---|---|---|---|---|
| $n$ | $PTM^{0.1}$ | SD | PM | HM | Mean |
| 20 | 19.688 | 17.746 | 14.392 | 13.851 | 1.0000 |
| | (121.34) | (134.61) | (165.99) | (172.47) | (2388.9) |
| 40 | 37.455 | 29.716 | 21.775 | 21.958 | 1.0000 |
| | (58.615) | (73.878) | (100.82) | (99.980) | (2195.4) |
| 60 | 54.039 | 39.880 | 27.763 | 27.848 | 1.0000 |
| | (39.687) | (53.778) | (77.249) | (77.014) | (2144.7) |
| 80 | 71.694 | 49.216 | 32.799 | 32.500 | 1.0000 |
| | (29.620) | (43.149) | (64.746) | (65.342) | (2123.6) |
| 100 | 83.543 | 53.948 | 35.536 | 35.198 | 1.0000 |
| | (24.760) | (38.343) | (58.210) | (58.770) | (2068.6) |

TABLE 4
*Finite sample efficiency of $PTM^\alpha$ relative to the sample mean*

| | $0.80N_2((0,0)', I) + 0.20N_2((10,10)', 25I)$ | | | | |
|---|---|---|---|---|---|
| $n$ | $PTM^{0.1}$ | SD | PM | HM | Mean |
| 20 | 51.779 | 37.573 | 28.653 | 27.416 | 1.0000 |
|    | (167.26) | (230.50) | (302.25) | (315.89) | (8660.5) |
| 40 | 89.745 | 52.678 | 35.790 | 34.733 | 1.0000 |
|    | (93.209) | (158.80) | (233.73) | (240.84) | (8356.0) |
| 60 | 121.86 | 62.637 | 40.018 | 39.876 | 1.0000 |
|    | (68.058) | (132.41) | (207.24) | (207.98) | (8293.4) |
| 80 | 155.80 | 68.658 | 43.001 | 42.364 | 1.0000 |
|    | (52.973) | (120.21) | (191.93) | (194.82) | (8253.1) |
| 100 | 176.39 | 71.282 | 43.807 | 43.385 | 1.0000 |
|     | (46.144) | (114.18) | (185.80) | (187.60) | (8139.3) |

**6. Concluding remarks.** We now account for the motivation of selecting projection depth for multivariate trimming, review some related trimmed means and studies in the literature, address the computing issues, discuss some practical choices of $\alpha$ values and summarize the major results obtained in this paper.

6.1. *Why projection depth trimmed means?* There are a number of depth notions in the literature; see [16], for example. For any given notion, one can define and study the corresponding depth trimmed means. Among the existing notions, the projection depth represents a very favorable one; see [33, 34, 37]. Tukey halfspace depth, also built based on projection pursuit methodology, is its major competitor. The projection depth, as a center-outward *strictly monotone* function, conveys more information about data points than the halfspace depth, a center-outward *step* function, does. As a matter of fact, the projection depth and its induced estimators can outperform the halfspace depth and its induced estimators with respect to two central performance criteria: robustness and efficiency. For example, the halfspace depth itself is much less robust than the projection depth (with appropriate choices of $\mu$ and $\sigma$). Indeed, the former has the lowest breakdown point ($\approx 0$) whereas the latter can have the highest breakdown point ($\approx 1/2$); see [34]. The estimators induced from the halfspace depth are also less robust than those from the projection depth. For example, the breakdown point of the halfspace median is $\leq 1/3$ [6] whereas that of the projection median can be about $1/2$ [33], the highest among all affine equivariant location estimators in high dimensions. The breakdown point of the $\alpha$th trimmed mean based on the halfspace depth is $(\lfloor \alpha n \rfloor + 1)/n$ whereas that of the one based on the projection depth is $\lfloor (n-d+1)/2 \rfloor/n$ (Theorem 3). On the other hand, the efficiency of the bivariate halfspace median relative to the mean is about 77% whereas that of the projection median can be as high as 95%; see [33]. The projection depth trimmed means can also be more efficient than those based on the halfspace depth. This is especially true when the underlying model slightly deviates from the assumed symmetric one (such as in the contaminated normal model cases). The robustness and efficiency advantages motivate us to focus on projection depth trimming. Note that with general choices of $(\mu, \sigma)$ we deal with a class of depth trimmed means instead of a single one as in the halfspace depth case. This is yet another motivation for projection depth trimming. The approaches and techniques in this paper, however, are applicable to other depth trimmed means and *covariance matrices* as well.

6.2. *Related estimators and studies in the literature.* First we note that by combining the integral regions with the integrands, (5) can be *trivially* written as

$$(12) \quad PTM^\alpha(F) = \int w^*(PD(x,F)) x \, dF(x) \Big/ \int w^*(PD(x,F)) \, dF(x),$$

with $w^*(s) = w(s) I(s \geq \alpha)$. Indeed we use this form repeatedly in the proofs. We adopt (5) [not (12)] since it is consistent with the regular univariate trimmed mean definition and manifests the depth trimming idea more clearly. Depth trimmed means with the form (12) have been discussed by Dümbgen [8] for simplicial depth, by Massé [18] for halfspace depth and by Zuo, Cui and He [35] for general depth. These discussions, however, are based on the assumption that $w^*(s)$ is continuously differentiable, which straightforwardly *excludes* (12) with $w^*(s) = w(s) I(s \geq \alpha)$. The difference here between continuous differentiability on a *closed* interval and discontinuity, seemingly very minor since it is *understood* that one can approximate the discontinuous function by a sequence of continuous differentiable ones, turns out to be *crucial*. The immediate problem with the sequence approach is the unbounded derivatives of the approximating functions. Boundedness is essential in the treatments of Dümbgen [8], Massé [18] and Zuo, Cui and He [35]. To deal with (alleviate) the unboundedness effect one is (essentially) forced to construct a *random* sequence depending on the convergence rate of the process $\sqrt{n}(PD(x, F_n) - PD(x, F))$. This, however, seems infeasible, or the halfspace depth median would have been asymptotically normal, which is not true, as shown in [1].

Note that with a nonzero constant $w$, (5) admits a 0–1 trimming scheme, which is the one used in the *regular* (univariate) *trimming*. This, however, is *not* the case with Dümbgen [8], Massé [18] and Zuo, Cui and He [35], where a continuous differentiable $w^*$ is assumed. This is yet another difference between this and the other papers.

Now in one dimension with the same 0–1 trimming scheme, this paper (5), introduces a *new* type of trimmed mean that is different from the regular univariate trimmed mean (which corresponds to halfspace depth trimming) as well as the metrically trimmed mean (see [2, 13]). Indeed, the trimming in this paper is based on the "generalized standardized deviation" (the outlyingness), whereas the regular trimming is based on the ranks of sample points. The metrically trimmed mean uses deviations to trim. But as in the regular trimming case, it always trims a fixed fraction of sample points. The projection depth trimming in this paper trims sample points only when they are "bad." The advantages of this trimming scheme include the gain in robustness (see comments after Theorem 3) and in efficiency for models which slightly deviate from the assumed symmetric ones.

Based on halfspace depth, Donoho and Gasko [5] introduced a trimmed mean [corresponding to (5) with a nonzero constant $w$ and $dF(x)$ replaced by $dx$] and studied its breakdown point; Nolan [19] and van der Vaart and Wellner [32] studied the asymptotic normality and the Hadamard-differentiability, respectively, of the *same* estimator. When introducing the notion of simplicial depth, Liu [14] also defined a depth trimmed mean, which is not based on depth regions, though.

6.3. *Computing projection depth trimmed means.* Like all other affine equivariant high-breakdown procedures, the projection depth trimmed means are computationally intensive. Exact computing in high dimensions, though possible (and an algorithm for two-dimensional data exists), is infeasible. Approximate computing is much faster and quite reliable and is sufficient in most applications. Basic approaches include randomly selecting projection directions or selecting among those perpendicular to the hyperplanes containing $d$ data points. Feasible approximate algorithms for high-dimensional data exist and are utilized in this paper.

6.4. *Choice of $\alpha$ values.* A very legitimate practical concern for $PTM^\alpha$ is the choice of the $\alpha$ value. Empirical evidence indicates that an $\alpha$ value around 0.05 to 0.1 can lead to a very efficient $PTM^\alpha$ at both light- and heavy-tailed distributions. One, of course, might also adopt an adaptive data-driven approach to determine an appropriate $\alpha$ value. For a given data set, an $\alpha$ value is determined based on the behavior of the tail of the data set. Generally speaking, a small value of $\alpha$ (e.g., 0.05) is selected for light-tailed data while a larger value is selected for a heavy-tailed one.

6.5. *Main results obtained in the paper.* This paper introduces projection depth trimmed means, examines their performance with respect to two principal criteria, robustness and efficiency, and establishes their limiting distribution via asymptotic representations. It turns out that the depth trimmed means can be highly robust locally (with bounded influence functions) as well as globally (with the best breakdown point among affine equivariant competitors). Robustness and efficiency do not work in tandem in general. Results obtained in the paper indicate, however, that the depth trimmed means, unlike the mean and the (halfspace) depth median, can keep a very good balance between the two. At normal and other light-tailed symmetric models, they (with high relative efficiency for suitable $\alpha$'s) are *better* choices than the depth medians and *strong* competitors to the sample mean, which is the best (in terms of efficiency) at the normal model. At contaminated (normal or symmetric) models (the more realistic ones in practice) and heavy-tailed models, they, with very robust and overwhelmingly high efficiency, are much better choices than the sample mean and other depth competitors. As a by-product, this paper introduces a new type of trimmed mean in one dimension which can have advantages over the regular and metrically trimmed means with respect to the two central performance criteria, robustness and efficiency.

## APPENDIX: SELECTED PROOFS AND AUXILIARY LEMMAS

PROOF OF THEOREM 1. To prove the theorem, we need the following lemmas.

LEMMA A.1. *Under* (C0) *and* (C2) *and for fixed* $x \in \mathbb{R}^d$ *and very small* $\varepsilon > 0$:

(a) $PD(y, F)$ *and* $PD(y, F(\varepsilon, \delta_x))$ *are Lipschitz continuous in* $y \in \mathbb{R}^d$;
(b) $\sup_{y \in \mathbb{R}^d}(1 + \|y\|)|PD(y, F(\varepsilon, \delta_x)) - PD(y, F)| = o_x(1)$;
(c) $\partial PD^\alpha(G) = \{y : PD(y, G) = \alpha\}$ *for* $0 < \alpha < \alpha^*$ *and* $G = F$ *or* $F(\varepsilon, \delta_x)$;
(d) $PD^{\alpha+\eta}(F) \subseteq PD^{\alpha+\eta/2}(F(\varepsilon, \delta_x)) \subseteq PD^\alpha(F)$ *for any* $0 < \eta \leq \alpha^* - \alpha$.

PROOF. Conditions (C0) and (C2) imply that $\sup_{\|u\|=1} |\mu(F_{\varepsilon u})|$ and $\sup_{\|u\|=1} \sigma(F_{\varepsilon u})$ are finite for sufficiently small $\varepsilon > 0$. It is readily seen that $|\inf_{\|u\|=1} \sigma(F_{\varepsilon u}) - \inf_{\|u\|=1} \sigma(F_u)| \leq \sup_{\|u\|=1} |\sigma(F_{\varepsilon u}) - \sigma(F_u)| = o_x(1)$. This, together with (C0), implies that $\inf_{\|u\|} \sigma(F_{\varepsilon u})$ is bounded below from 0 for fixed $x$ and sufficiently small $\varepsilon$. Hence for sufficiently small $\varepsilon > 0$,

(C0′) $\sup_{\|u\|=1} \mu(F_{\varepsilon u}) < \infty, \ 0 < \inf_{\|u\|} \sigma(F_{\varepsilon u}) \leq \sup_{\|u\|=1} \sigma(F_{\varepsilon u}) < \infty$.

The Lipschitz continuity of $PD(\cdot, F)$ can be established by following the proof of Theorem 2.2 of [33]. For that of $PD(\cdot, F_\varepsilon)$, we observe that for small $\varepsilon > 0$

$$|PD(y_1, F_\varepsilon) - PD(y_2, F_\varepsilon)| \leq |O(y_1, F_\varepsilon) - O(y_2, F_\varepsilon)|$$
$$\leq \|y_1 - y_2\| / \inf_{\|u\|=1} \sigma(F_{\varepsilon u}).$$

This and (C0′) lead to the Lipschitz continuity of $PD(\cdot, F_\varepsilon)$. Part (a) follows.

To show part (b), first we observe that

$$(1 + \|y\|)|PD(y, F_\varepsilon) - PD(y, F)|$$
$$< \frac{1 + \|y\|}{1 + O(y, F_\varepsilon)}$$
$$\times \frac{O(y, F)\sup_{\|u\|=1}|\sigma(F_{\varepsilon u}) - \sigma(F_u)| + \sup_{\|u\|=1}|\mu(F_{\varepsilon u}) - \mu(F_u)|}{\inf_{\|u\|=1}\sigma(F_{\varepsilon u})(1 + O(y, F))}$$
$$\leq \frac{1 + \|y\|}{1 + |\|y\| - \mu(F_{\varepsilon y_0})|/(\sup_{\|u\|}\sigma(F_{\varepsilon u}))}$$
$$\times \frac{\sup_{\|u\|=1}|\sigma(F_{\varepsilon u}) - \sigma(F_u)| + \sup_{\|u\|=1}|\mu(F_{\varepsilon u}) - \mu(F_u)|}{\inf_{\|u\|=1}\sigma(F_{\varepsilon u})},$$

where $y_0 = y/\|y\|$. Part (b) now follows immediately from (C0′) and (C2).

Part (c) with $G = F$ is covered by Theorem 2.3 of [33]. To show the case $G = F_\varepsilon$, first we note that $PD^\alpha(F_\varepsilon)$ is nonempty for sufficiently small $\varepsilon > 0$ since by (b) for any $\theta \in \mathbb{R}^d$ with $PD(\theta, F) = \alpha^*$, $PD(\theta, F_\varepsilon) > \alpha$ for small $\varepsilon > 0$.

We now show that $\{y : PD(y, F_\varepsilon) = \alpha\} \subseteq \partial PD^\alpha(F_\varepsilon)$, the boundary of $PD^\alpha(F_\varepsilon)$. Let $PD(y, F_\varepsilon) = \alpha$. Such $y$ exists since (i) by (C0′) we can show that $PD(z, F_\varepsilon) \to 0$ as $\|z\| \to \infty$ (see Theorem 2.1 of [33]) and (ii) $PD(\cdot, F_\varepsilon)$ is Lipschitz continuous by (a) and $PD(\theta, F_\varepsilon) > \alpha$ for sufficiently small $\varepsilon > 0$. Assume that $y \notin \partial PD^\alpha(F_\varepsilon)$; that is, $y$ is an interior point of $PD^\alpha(F_\varepsilon)$. Then there is a small ball centered at $y$ with radius $r$ and contained in the interior of $PD^\alpha(F_\varepsilon)$. By the scale equivariance of $\mu$ we see immediately that there is a direction $u_0$ such that

$$(u_0'y - \mu(F_{\varepsilon u_0}))/\sigma(F_{\varepsilon u_0}) > O(y, F_\varepsilon) - r/\sup_{\|u\|=1}\sigma(F_{\varepsilon u})$$

for sufficiently small $\varepsilon$ such that $\sup_{\|u\|=1}\sigma(F_{\varepsilon u}) < \infty$. On the other hand, we see that $y' = y + u_0 r \in PD^\alpha(F_\varepsilon)$ and

$$O(y', F_\varepsilon) \geq \frac{u_0'y' - \mu(F_{\varepsilon u_0})}{\sigma(F_{\varepsilon u_0})} = \frac{u_0'y - \mu(F_{\varepsilon u_0})}{\sigma(F_{\varepsilon u_0})} + \frac{r}{\sigma(F_{\varepsilon u_0})} > O(y, F_\varepsilon).$$

But this implies that $PD(y', F_\varepsilon) < PD(y, F_\varepsilon) = \alpha$, which is a contradiction.

We now show that $\partial PD^\alpha(F_\varepsilon) \subseteq \{y : PD(y, F_\varepsilon) = \alpha\}$. Let $y \in \partial PD^\alpha(F_\varepsilon)$. Then by the continuity of $PD(\cdot, F_\varepsilon)$ for sufficiently small $\varepsilon > 0$, we conclude that $PD(y, F_\varepsilon) \geq \alpha$. If $PD(y, F_\varepsilon) > \alpha$, then by the continuity of $PD(\cdot, F_\varepsilon)$ there is a small ball $B(y)$ centered at $y$ with $PD(z, F_\varepsilon) > \alpha$ for all $z \in B(y)$ for sufficiently small $\varepsilon > 0$. But this contradicts the assumption that $y \in \partial PD^\alpha(F_\varepsilon)$. Part (c) follows.

By part (a), for the given $0 < \eta \le \alpha^* - \alpha$ and any $y \in PD^{\alpha+\eta/2}(F_\varepsilon)$

$$PD(y, F) \ge PD(y, F_\varepsilon) - \eta/2 \ge \alpha + \eta/2 - \eta/2 = \alpha$$

for sufficiently small $\varepsilon > 0$. Likewise, for any $y \in PD^{\alpha+\eta}(F)$, by part (a)

$$PD(y, F_\varepsilon) \ge PD(y, F) - \eta/2 \ge \alpha + \eta - \eta/2 = \alpha + \eta/2$$

for sufficiently small $\varepsilon > 0$. Part (d) follows immediately. $\square$

LEMMA A.2. *Under* (C1)–(C2), $\sup_{\|u\|=1} |R^\alpha(u, F(\varepsilon, \delta_x)) - R^\alpha(u, F)| = o_x(1)$.

PROOF. By (C1), for any $x \in \mathbb{R}^d$ there is a unit vector $v(x)$ such that $g(x, v(x), F) = O(x, F)$, where $g(x, \cdot, F)$ is defined before (2). Let $v(u) := v(R(u)u)$. By (C1) [hence (C0)] and Lemma A.1, $O(R(u)u, F) = \beta(\alpha)$. Thus

$$R(u)u'v(u) = \beta(\alpha)\sigma(F_{v(u)}) + \mu(F_{v(u)}), \qquad R(u)u'v \le \beta(\alpha)\sigma(F_v) + \mu(F_v)$$

for any $v \in S^{d-1}$. Likewise, for $v_\varepsilon(u) := v(R^\alpha(u, F(\varepsilon, \delta_x))u)$ and small $\varepsilon > 0$

$$g(R^\alpha(u, F(\varepsilon, \delta_x))u, v_\varepsilon(u), F(\varepsilon, \delta_x)) = O(R^\alpha(u, F(\varepsilon, \delta_x))u, F(\varepsilon, \delta_x)) = \beta(\alpha).$$

Again for convenience, write $R_\varepsilon(u)$ or $R_\varepsilon^\alpha(u)$ for $R^\alpha(u, F(\varepsilon, \delta_x))$. Hence we have

$$R_\varepsilon(u)u'v_\varepsilon(u) = \beta(\alpha)\sigma(F_{\varepsilon v_\varepsilon(u)}) + \mu(F_{\varepsilon v_\varepsilon(u)}),$$

$$R_\varepsilon(u)u'v \le \beta(\alpha)\sigma(F_{\varepsilon v}) + \mu(F_{\varepsilon v})$$

for any unit vector $v \in S^{d-1}$. These and the counterparts above yield

$$(\beta(\alpha)[\sigma(F_{\varepsilon v_\varepsilon(u)}) - \sigma(F_{v_\varepsilon(u)})] + (\mu(F_{\varepsilon v_\varepsilon(u)}) - \mu(F_{v_\varepsilon(u)})))/u'v_\varepsilon(u)$$
$$\le R_\varepsilon(u) - R(u)$$
$$\le (\beta(\alpha)[\sigma(F_{\varepsilon v(u)}) - \sigma(F_{v(u)})] + (\mu(F_{\varepsilon v(u)}) - \mu(F_{v(u)})))/u'v(u).$$

If we can show that both $\inf_{\|u\|} |u'v(u)|$ and $\inf_{\|u\|=1} |u'v_\varepsilon(u)|$ are bounded away from 0, the desired result follows immediately from (C2).

Since $PD^\alpha(F)$ is assumed to contain the origin, the deepest point, thus $O(0, F) = \beta(\alpha^*) < \beta(\alpha) = O(R(u)u, F)$ for any $\|u\| = 1$ (Theorem 2.3 of [33]). Hence

$$\beta(\alpha) = \frac{u'v(u)R(u) - \mu(F_{v(u)})}{\sigma(F_{v(u)})} > O(0, F) + \gamma(\alpha, \alpha^*)$$

$$= \sup_{\|u\|=1} \frac{\mu(F_u)}{\sigma(F_u)} + \gamma(\alpha, \alpha^*)$$

for $\gamma(\alpha, \alpha^*) := (\beta(\alpha) - \beta(\alpha^*))/2$ and any $u \in S^{d-1}$, which, in turn, implies that

$$\frac{u'v(u)R(u)}{\sigma(F_{v(u)})} > \sup_{\|u\|=1} \frac{\mu(F_u)}{\sigma(F_u)} + \frac{\mu(F_{v(u)})}{\sigma(F_{v(u)})} + \gamma(\alpha, \alpha^*) \geq \gamma(\alpha, \alpha^*).$$

Note that $PD^\alpha(F)$ is bounded with a nonempty interior (see Theorem 2.3 of [33]). Hence $0 < R(u) < \infty$ uniformly in $u$. Now we have

$$\inf_{\|u\|=1} |u'v(u)| \geq \left\{ \inf_{\|u\|} \sigma(F_u)(\gamma(\alpha, \alpha^*)) \right\} \Big/ \sup_{\|u\|=1} R(u) > 0.$$

The argument for showing $\inf_{\|u\|=1} |u'v_\varepsilon(u)| > 0$ is the same. First we have

$$0 \in PD^{(\alpha+\delta)}(F) \subseteq PD^{(\alpha+\delta/2)}(F_\varepsilon) \subseteq PD^\alpha(F_\varepsilon),$$

by Lemma A.1 for some $0 < \delta < \alpha^* - \alpha$ and sufficiently small $\varepsilon$, which gives

$$O(0, F_\varepsilon) < O(R^{(\alpha+\delta/2)}(u, F_\varepsilon)u, F_\varepsilon) < O(R^\alpha(u, F_\varepsilon)u, F_\varepsilon) = \beta(\alpha),$$

uniformly in the unit vector $u$ for sufficiently small $\varepsilon > 0$ by (c) of Lemma A.1. Now treating $O(0, F_\varepsilon)$ as the $\beta(\alpha^*)$ above, we have the desired result by virtue of (C1)–(C2), (C0$'$), Lemma A.1 and Theorem 2.3 of [33]. □

LEMMA A.3. (a) $R^\alpha(u, F)$ *is continuous in* $u$ *if* (C0) *holds.* (b) $R^\alpha(u, F_\varepsilon)$ *is continuous in* $u$ *for sufficiently small* $\varepsilon > 0$ *if* (C0) *and* (C2) *hold.*

PROOF. (a) Suppose that $R^\alpha(u, F)$ is not continuous in $u$. Then there is a sequence $u_m \to u_0$ such that $\limsup_{m \to \infty} R^\alpha(u_m, F) \neq R^\alpha(u_0, F)$. By the boundedness of $PD^\alpha(F)$, there is a subsequence $u_{m_k}$ of $u_m$ such that $u_{m_k} \to u_0$ and $\lim_{k \to \infty} R^\alpha(u_{m_k}, F) = R_0^\alpha \neq R^\alpha(u_0, F)$. Note that $R_0^\alpha$ must be less than $R^\alpha(u_0, F)$ since otherwise we have by the uniform continuity of $PD(x, F)$ in $x$ that

$$\lim_{k \to \infty} PD(R^\alpha(u_{m_k}, F)u_{m_k}, F) = PD(R_0^\alpha u_0, F) = \alpha = PD(R^\alpha(u_0, F)u_0, F),$$

which contradicts the definition of $R^\alpha(u_0, F)$. Thus $R_0^\alpha < R^\alpha(u_0, F)$. The quasi-concavity of $PD(\cdot, F)$ (Theorem 2.1 of [33]) implies $PD(x, F) = \alpha$ for any point $x \in [R_0^\alpha u_0, R^\alpha(u_0, F)u_0]$, and further all such points $x$ are boundary points of $PD^\alpha(F)$ in light of Lemma A.1.

Let $N_0(\varepsilon_0)$ be a small ball centered at the deepest point, the origin, and contained in $PD^\alpha(F)$. Let $N_{x_0}(\varepsilon_1)$ be a small ball centered at $x_0 \in [R_0^\alpha u_0, R^\alpha(u_0, F)u_0]$ and $\varepsilon_1$ small enough such that the ray stemming from $R^\alpha(u_0, F)u_0$ and passing through the ball $N_{x_0}(\varepsilon_1)$ always passes through the ball $N_0(\varepsilon_0)$. But then there is a point $y_0 \in N_{x_0}(\varepsilon_1)$ and $y_0 \notin PD^\alpha(F)$ and a point $y_1 \in N_0(\varepsilon_0)$ such that

$$y_0 = \lambda y_1 + (1 - \lambda) R^\alpha(u_0, F) u_0, \qquad PD(y_0, F) < \alpha,$$

for some $0 < \lambda < 1$, which contradicts the quasi-concavity of $PD(\cdot, F)$.

(b) From the proof of Lemma A.1 we see that (C0$'$) holds for sufficiently small $\varepsilon > 0$ by virtue of (C0) and (C2). Then the quasi-concavity of $PD(\cdot, F_\varepsilon)$ follows:

$$PD(\lambda x + (1-\lambda)y, F_\varepsilon) \geq \min\{PD(x, F_\varepsilon) PD(y, F_\varepsilon)\}$$

for any $0 < \lambda < 1$ and sufficiently small $\varepsilon > 0$. Now invoking Lemma A.1 and the arguments utilized in (a) we can complete the proof. □

We now prove Theorem 1. Following the proof of Lemma A.2, we have

$$(\beta(\alpha)[\sigma(F_{\varepsilon v_\varepsilon(u)}) - \sigma(F_{v_\varepsilon(u)})] + (\mu(F_{\varepsilon v_\varepsilon(u)}) - \mu(F_{v_\varepsilon(u)})))/u'v_\varepsilon(u)$$
$$\leq R_\varepsilon(u) - R(u)$$
$$\leq (\beta(\alpha)[\sigma(F_{\varepsilon v(u)}) - \sigma(F_{v(u)})] + (\mu(F_{\varepsilon v(u)}) - \mu(F_{v(u)})))/u'v(u),$$

and $\inf_{\|u\|=1} |u'v(u)|$ and $\inf_{\|u\|=1} |u'v_\varepsilon(u)|$ are bounded below from 0 for sufficiently small $\varepsilon > 0$. The desired result then follows from (C3) and the continuity of $IF(v(u)'x; \sigma, F_{v(u)})$ and $IF(v(u)'x; \mu, F_{v(u)})$ in $v(u)$ for $u \in A$, provided that we can show further that $v_\varepsilon(u) \to v(u)$ uniformly in $u$ as $\varepsilon \to 0$.

We first show that $v_\varepsilon(u) \to v(u)$ as $\varepsilon \to 0$ for a fixed $u$. If it is not true, then there are a sequence $\varepsilon_n \to 0$ and a small $\eta > 0$ such that $\|v_{\varepsilon_n}(u) - v(u)\| \geq \eta$ for $n \geq 1$. By the compactness of $S^{d-1}$, there is a subsequence of $v_{\varepsilon_n}(u)$, denoted still by $v_{\varepsilon_n}(u)$ for simplicity, that converges to $v_0 \in S^{d-1}$. Observe that

(13) $O(R^\alpha(u, F_{\varepsilon_n})u, F_{\varepsilon_n}) = (u'v_{\varepsilon_n}(u)R^\alpha(u, F_{\varepsilon_n}) - \mu(F_{\varepsilon_n v_{\varepsilon_n}(u)}))/\sigma(F_{\varepsilon_n v_{\varepsilon_n}(u)}).$

Following the proof of Theorem 2.2 of [33] and by (C1) and (C2), we have (i) the Lipschitz continuity of $O(\cdot, F_{\varepsilon_n})$ for small $\varepsilon > 0$ and (ii) for large $M > 0$

(14) $$\sup_{\|y\| < M} |O(y, F_{\varepsilon_n}) - O(y, F)| \to 0 \quad \text{as } n \to \infty.$$

These, together with (13), Lemmas A.1 and A.2 and (C1)–(C2), yield

$$\frac{u'v_0 R^\alpha(u, F) - \mu(F_{v_0})}{\sigma(F_{v_0})} = O(R^\alpha(u, F)u, F) = \frac{u'v(u) R^\alpha(u, F) - \mu(F_{v(u)})}{\sigma(F_{v(u)})}.$$

Uniqueness of $v(u) = v(y)$ for $y = R^\alpha(u, F)u$ implies that $v(u) = v_0$, which, however, contradicts $\|v_0 - v(u)\| \geq \eta$. Hence $v_\varepsilon(u) \to v(u)$ for any fixed $u \in S^{d-1}$.

With the same argument, we can show that the convergence is uniform in the unit vector $u$, since otherwise there are a sequence $u_n \in S^{d-1}$, a sequence $\varepsilon_n$ ($\varepsilon_n \downarrow 0$ as $n \to \infty$) and some small $\eta > 0$ such that $|v_{\varepsilon_n}(u_n) - v(u_n)| > \eta$ as $n \to \infty$. By the compactness of $S^{d-1}$, there is a subsequence $u_{n_m}$ of $u_n$ such that $u_{n_m} \to u_0 \in S^{d-1}$, $v(u_{n_m}) \to v_0 \in S^{d-1}$ and $v_{\varepsilon_{n_m}}(u_{n_m}) \to v_1 \in S^{d-1}$ as $m \to \infty$. We then can show that $v_0 = v_1 = v(u_0)$ (here we need Lemma A.3). But this contradicts $|v_{\varepsilon_n}(u_{n_m}) - v(u_{n_m})| > \eta$ as $m \to \infty$. The desired result follows. □

PROOF OF COROLLARY 1. We first verify the conditions in Theorem 1. We see that $\mu(F_u)(=0)$ and $\sigma(F_u)(=\sqrt{u'\Sigma u}m_0 > 0)$ are continuous in $u \in S^{d-1}$ and

$$\mu(F_u(\varepsilon, \delta_x)) = \text{Med}\{\sqrt{u'\Sigma u}F_Z^{-1}(a_\varepsilon), u'x, \sqrt{u'\Sigma u}F_Z^{-1}(b_\varepsilon)\},$$

$$\sigma(F_u(\varepsilon, \delta_x)) = \text{Med}\{\sqrt{u'\Sigma u}F_{Z_\varepsilon}^{-1}(a_\varepsilon), |u'x - \mu(F_{\varepsilon u})|, \sqrt{u'\Sigma u}F_{Z_\varepsilon}^{-1}(b_\varepsilon)\},$$

where $a_\varepsilon = (1-2\varepsilon)/(2(1-\varepsilon))$, $b_\varepsilon = 1/(2(1-\varepsilon))$ and $Z_\varepsilon = |Z - \mu(F_{\varepsilon u})/\sqrt{u'\Sigma u}|$. It follows that both $\mu(F_u(\varepsilon, \delta_x))$ and $\sigma(F_u(\varepsilon, \delta_x))$ are continuous in $u \in S^{d-1}$ for sufficiently small $\varepsilon > 0$. Thus (C1) holds. The last two displays also lead to (C2):

$$\sup_{\|u\|=1} |\mu(F_u(\varepsilon, \delta_x)) - \mu(F_u)| = o_x(1), \qquad \sup_{\|u\|=1} |\sigma(F_u(\varepsilon, \delta_x)) - \sigma(F_u)| = o_x(1).$$

Note that for any $y \in \partial PD^\alpha(F)$ it can be seen that $v(y) = \Sigma^{-1}y/\|\Sigma^{-1}y\|$, $O(y, F) = \|\Sigma^{-1/2}y\|/m_0 = \beta(\alpha)$ and $PD(y, F) = (1 + O(y, F))^{-1} = \alpha$. Hence $\mathcal{U}(y) = \{v(y) : g(y, v(y), F) = O(y, F)\}$ is a singleton for any $y \in \partial PD^\alpha(F)$.

It can be shown that for any $u \in S^{d-1}$

(15)   $IF(u'x; \mu, F_u) = \sqrt{u'\Sigma u}\,\text{sign}(u'x)/(2h_z(0))$,

(16)   $IF(u'x; \sigma, F_u) = \sqrt{u'\Sigma u}\,\text{sign}(|u'x| - \sqrt{u'\Sigma u}m_0)/(4h_z(m_0))$.

These and the expressions for $\mu(F_u(\varepsilon, \delta_x))$ and $\sigma(F_u(\varepsilon, \delta_x))$ above lead to (C3).

Obviously, both $IF(u'x; \mu, F_u)$ and $IF(u'x; \sigma, F_u)$ are continuous in $u \in S^{d-1}$ if $x = 0$. When $x \neq 0$, $IF(v(y)'x; \mu, F_{v(y)})$ is continuous in $v(y)$ for any $y \in A^* \subseteq \partial PD^\alpha(F)$ with $\partial PD^\alpha(F) - A^* = \{y : y'\Sigma^{-1}x = 0, y \in \partial PD^\alpha(F)\}$ and $P(\{y : y'\Sigma^{-1}x = 0, y \in \partial PD^\alpha(F)\}) = 0$ for fixed $x \in \mathbb{R}^d$. Likewise, when $x \neq 0$ we see that $IF(v(y)'x; \sigma, F_{v(y)})$ is continuous in $v(y)$ for any $y \in A^{**} \subseteq \partial PD^\alpha(F)$ with $\partial PD^\alpha(F) - A^{**} = \{y : y'\Sigma^{-1}x = \pm\beta(\alpha)m_0^2, y \in \partial PD^\alpha(F)\}$. The latter set is empty if $\|\Sigma^{-1/2}x\| < m_0$. Also $P(\{y : y'\Sigma^{-1}x = \pm\beta(\alpha)m_0^2, y \in \partial PD^\alpha(F)\}) = 0$ for fixed $x \in \mathbb{R}^d$. Thus there is a set $A \subseteq S^{d-1}$ with $P\{y : v(y) \in S^{d-1} - A, y \in \partial PD^\alpha(F)\} = 0$ such that $IF(v(u)'x; \mu, F_{v(u)})$ and $IF(v(u)'x; \sigma, F_{v(u)})$ with $v(u) = v(R^\alpha(u, F)u)$ are continuous in $v(u)$ for all $u \in A$. Here $A = S^{d-1}$ if $x = 0$ and $A = S^{d-1} - \{u : u'\Sigma^{-1}x = 0\} \cup \{u : u'\Sigma^{-1}x = \pm\|\Sigma^{-1/2}u\|m_0\}$ if $x \neq 0$.

Invoking Theorem 1 and (15) and (16), we have the desired result. □

PROOF OF THEOREM 2. To prove the theorem, we need the following lemma.

LEMMA A.4. *Under* (C1) *and* (C3)–(C5), *we have*

(17) $$(PD(y, F(\varepsilon, \delta_x)) - PD(y, F))/\varepsilon = h(x, y) + o_x(1) \quad \text{uniformly in } y \in B.$$

PROOF. By the conditions and the proof of Lemma 6.1 of [36], we have for $u(y, \tau) = \{u : \|u\| = 1, \|u - v(y)\| \leq \tau\}$, $\tau > 0$,

$$\inf_{u \in u(y,\tau)} (-g(y, u, F_\varepsilon) + g(y, u, F))$$

$$\leq -O(y, F_\varepsilon) + O(y, F) \leq -g(y, v(y), F_\varepsilon) + g(y, v(y), F),$$

for $y \in B$. The given conditions on $\mu$ and $\sigma$ imply (C0). Hence $PD^{\alpha-\eta}(F)$ is bounded (Theorem 2.3 of [33]). Conditions (C3) and (C5) imply that $\sigma(F_{\varepsilon u}) \to \sigma(F_u)$ uniformly for $u \in \{v(y) : y \in B\}$. These and (C3) yield

$$\frac{-g(y, u, F_\varepsilon) + g(y, u, F)}{\varepsilon(1 + O(y, F))}$$

$$= \frac{g(y, u, F) IF(u'x; \sigma, F_u) + IF(u'x; \mu, F_u)}{\sigma(F_u)(1 + O(y, F))} + o_x(1)$$

uniformly in $y \in \mathbb{R}^d$ and in $u \in S^{d-1}$. Hence we have

$$\inf_{u \in u(y,\tau)} \frac{g(y, u, F) IF(u'x; \sigma, F_u) + IF(u'x; \mu, F_u)}{\sigma(F_u)(1 + O(y, F))(1 + O(y, F_\varepsilon))} + o_x(1)$$

$$\leq \frac{PD(y, F_\varepsilon) - PD(y, F)}{\varepsilon}$$

$$\leq \frac{g(y, v(y), F) IF(v(y)'x; \sigma, F_{v(y)}) + IF(v(y)'x; \mu, F_{v(y)})}{\sigma(F_{v(y)})(1 + O(y, F))(1 + O(y, F_\varepsilon))} + o_x(1)$$

uniformly in $y$ over $B$. Let $\tau = \varepsilon/2$. By the given conditions, the result follows. □

We now prove Theorem 2 based on the lemma. First we can write for fixed $\alpha$,

(18) $$PTM(F_\varepsilon) - PTM(F) = \frac{\int_{PD^\alpha(F_\varepsilon)} (y - PTM(F)) w(PD(y, F_\varepsilon)) \, dF_\varepsilon(y)}{\int_{PD^\alpha(F_\varepsilon)} w(PD(y, F_\varepsilon)) \, dF_\varepsilon(y)}.$$

The denominator can be written as

$$(1 - \varepsilon) \int I(y \in PD^\alpha(F_\varepsilon)) w(PD(y, F_\varepsilon)) \, dF(y)$$

$$+ \varepsilon I(PD(x, F_\varepsilon) \geq \alpha) w(PD(x, F_\varepsilon)),$$

which, with Lemma A.1 and Lebesgue's dominated convergence theorem, yields

$$(19) \quad \int_{PD^\alpha(F_\varepsilon)} w(PD(y, F_\varepsilon)) \, dF_\varepsilon(y) = \int_{PD^\alpha(F)} w(PD(y, F)) \, dF(y) + o_x(1).$$

The numerator of (18) can be decomposed into three parts,

$$I_{1\varepsilon} = \int_{PD^\alpha(F_\varepsilon)} y^* w(PD(y, F_\varepsilon)) \, dF_\varepsilon(y) - \int_{PD^\alpha(F)} y^* w(PD(y, F_\varepsilon)) \, dF_\varepsilon(y),$$

$$I_{2\varepsilon} = \int_{PD^\alpha(F)} y^* w(PD(y, F_\varepsilon)) \, dF_\varepsilon(y) - \int_{PD^\alpha(F)} y^* w(PD(y, F)) \, dF_\varepsilon(y),$$

$$I_{3\varepsilon} = \int_{PD^\alpha(F)} y^* w(PD(y, F)) \, dF_\varepsilon(y),$$

where $y^* = y - PTM(F)$. It follows immediately that

$$(20) \quad I_{3\varepsilon}/\varepsilon = I(PD(x, F) \geq \alpha)(x - PTM(F))w(PD(x, F)).$$

(C1) implies (C0). This, the continuity of $w^{(1)}(\cdot)$ and Lemma A.1 yield

$$(21) \quad w(PD(y, F_\varepsilon)) - w(PD(y, F)) = (w^{(1)}(PD(y, F)) + o_x(1))H(y, F_\varepsilon),$$

uniformly in $y$ for the given $x$, where $H(y, F_\varepsilon) = PD(y, F_\varepsilon) - PD(y, F)$. By Lemma A.4, (C5) and the boundedness of $PD^\alpha(F)$, it is seen that $(1 + \|y\|)IF(x; PD(y, F), F)$ is bounded uniformly in $y$ for $y \in B$ and the given $x \in \mathbb{R}^d$. This, together with (21), Lemma A.4 and the boundedness of $PD^\alpha(F)$, immediately gives

$$(22) \quad I_{2\varepsilon}/\varepsilon = \int_{PD^\alpha(F)} (y - PTM(F))w^{(1)}(PD(y, F))h(x, y) \, dF(y) + o_x(1).$$

Write $\Delta(y, \varepsilon, \alpha)$ for $I(PD(y, F_\varepsilon) \geq \alpha)) - I(PD(y, F) \geq \alpha))$. By virtue of (21), Lemmas A.1 and A.4, the boundedness of $PD^\alpha(F)$ and $PD^\alpha(F_\varepsilon)$ for small $\varepsilon > 0$, and the argument used in the denominator of (18), we have

$$\frac{1}{\varepsilon} I_{1\varepsilon} = \frac{1}{\varepsilon} \int \Delta(y, \varepsilon, \alpha) y^* w(PD(y, F)) \, dF(y)$$

$$- \int \Delta(y, \varepsilon, \alpha) y^* w(PD(y, F)) \, dF(y)$$

$$+ \int \Delta(y, \varepsilon, \alpha) y^* w^{(1)}(PD(y, F)) \frac{PD(y, F_\varepsilon) - PD(y, F)}{\varepsilon} \, dF(y) + o_x(1).$$

Call the last three (integral) terms $I_{1\varepsilon i}$, $i = 1, 2, 3$, respectively. Then by Lemma A.1, (C3), (C5), the boundedness of $PD^\alpha(F)$ and $PD^\alpha(F_\varepsilon)$ for small $\varepsilon > 0$, the condition on $w$ and Lebesgue's dominated convergence theorem $I_{1\varepsilon 2} = o_x(1)$ and $I_{1\varepsilon 3} = o_x(1)$. For $I_{1\varepsilon 1}$, by the mean value theorem and Theorem 1 we have

$$I_{1\varepsilon 1} = \frac{1}{\varepsilon} \int_{S^{d-1}} \left[ \int_{R^\alpha(u,F)}^{R^\alpha(u,F_\varepsilon)} (ru - PTM(F)) w(PD(ru,F)) |J(u,r)| f(ru) \, dr \right] du$$

$$= \int_{S^{d-1}} ((\theta_\varepsilon(u) u - PTM(F)) w(PD(\theta_\varepsilon(u) u, F)) |J(u, \theta_\varepsilon(u))| f(\theta_\varepsilon(u) u)$$

$$\times (IF(x; R^\alpha(u,F), F) + o_x(1))) \, du,$$

where $\theta_\varepsilon(u)$ is a point between $R^\alpha(u, F_\varepsilon)$ and $R^\alpha(u, F)$ and $o_x(1)$ is in the uniform sense with respect to $u$. By Lemmas A.1 and A.2, (C5), the conditions on $f$ and $w$, the structure of $J(u,r)$, the boundedness of $PD^\alpha(F)$ and $PD^\alpha(F_\varepsilon)$ for small $\varepsilon > 0$ and Lebesgue's dominated convergence theorem, we have for $R^*(u) = R(u)u - PTM$,

$$\frac{1}{\varepsilon} I_{1\varepsilon} = \int_{S^{d-1}} R^*(u) w(PD(R(u)u, F)) |J(u, R(u))|$$

$$\times f(R(u)u) IF(x; R(u), F) \, du + o_x(1).$$

The desired result now follows immediately from this, (19), (20) and (22). □

PROOF OF THEOREM 3. Write $u'X^n$ for $F_{nu}$ and $X^n$ for $F_n$ and skip the $d=1$ case.

Consider the case $d > 1$. We first show $m = \lfloor (n-d+1)/2 \rfloor$ *contaminating points are enough to break down* $PTM^\alpha$. Move $m$ points of $X^n$ to the same site $y$. Denote the resulting data $X_m^n = \{Z_1, \ldots, Z_n\}$. Assume the first $m$ points $Z_i$ ($1 \le i \le m$) are at site $y$ far away from the cloud $X^n$. For $u \in S^{d-1}$, the projected data set (to direction $u$) is $\{u'Z_1, \ldots, u'Z_n\}$. Since $m + \lfloor (n+d+2)/2 \rfloor > n$, thus $|u'Z_i - \mu(u'X_m^n)|/\sigma(u'X_m^n) \le 2$ for all $1 \le i \le m$. This implies that $O(Z_i, X_m^n) \le 2$ for all $1 \le i \le m$. Hence $Z_i \in PD^\alpha(X_m^n)$ for $1 \le i \le m$ by (11). Since $\|\sum_{i=1}^m Z_i w(PD(Z_i, X_m^n))\| \to \infty$ as $\|y\| \to \infty$, therefore $PTM^\alpha(X_m^n)$ breaks down.

Now we show that $m = \lfloor (n-d+1)/2 \rfloor - 1$ *contaminating points are not enough to break down* $PTM^\alpha$. Again let $X_m^n = \{Z_1, \ldots, Z_n\}$ be any contaminated data set. Since $m < \lfloor (n+1)/2 \rfloor$ and $m + \lfloor (n+d+2)/2 \rfloor \le n$, $\sup_{\|u\|=1} \mu(u'X_m^n) < \infty$ and $\sup_{\|u\|=1} \sigma(u'X_m^n) < \infty$ uniformly for any contaminated data $X_m^n$ with $m$ original points contaminated. Hence $O(y, X_m^n) \to \infty$ as $\|y\| \to \infty$. That is, $y \notin PD^\alpha(X_m^n)$ when $\|y\|$ becomes very large. So $PTM^\alpha$ will not break down unless $PD^\alpha(X_m^n) \cap X_m^n$ becomes empty. We now show that the latter cannot happen.

Denote $\mu_u = \mu(u'X_m^n)$, $\sigma_u = \sigma(u'X_m^n)$ and $n_\sigma = \lfloor (n+d+2)/2 \rfloor$. Let $|u'Z_{i1} - \mu_u| \le \cdots \le |u'Z_{in_\sigma} - \mu_u| \le \cdots \le |u'Z_{in} - \mu_u|$ with the understanding that $\mu_u$, $\sigma_u$ and $Z_{ij}$ depend on $X_m^n$ and $u$ for all $1 \le j \le n$. Since $m + d + 1 \le n_\sigma$,

hence among $Z_{i1}, \ldots, Z_{in_\sigma}$ there are at least $d+1$ original points from $X^n$. Therefore

$$\sigma(u'X_m^n) \geq \inf_{i_1,\ldots,i_{d+1}} \max_{1 \leq k,l \leq (d+1)} |u'(X_{i_k} - X_{i_l})|/2$$

for any $X_m^n$ and $u$. Here $i_1, \ldots, i_r$ are $r$ arbitrary distinct integers from $\{1, \ldots, n\}$.

Clearly, there are at least $m_0 = \lfloor (n+2)/2 \rfloor + 1$ original points, say (without loss of generality) $X_i$, $1 \leq i \leq m_0$, uncontaminated. Then it is not difficult to see that

$$|u'X_i - \mu(u'X_m^n)| \leq \max_{i_1,\ldots,i_{(m_0-1)}} \max_{1 \leq k,l \leq (m_0-1)} |u'(X_{i_k} - X_{i_l})|, \qquad 1 \leq i \leq m_0,$$

for any $X_m^n$ and $u$. This, in conjunction with the last display, immediately yields

$$O(X_i; X_m^n) \leq \sup_{\|u\|=1} \frac{\max_{i_1,\ldots,i_{(m_0-1)}} \max_{1 \leq k,l \leq (m_0-1)} |u'(X_{i_k} - X_{i_l})|}{\inf_{i_1,\ldots,i_{d+1}} \max_{1 \leq k,l \leq (d+1)} |u'(X_{i_k} - X_{i_l})|/2},$$

for $1 \leq i \leq m_0$. Hence $X_i \in PD^\alpha(X_m^n)$ for all $1 \leq i \leq m_0$ for any $0 < \alpha < \alpha_d$. That is, $PD^\alpha(X_m^n) \cap X_m^n$ is not empty. We complete the proof. $\square$

PROOF OF THEOREM 4. The following lemma, an analogue of Lemma A.1, is needed in the sequel. It can be proved in much the same way as Lemma A.1.

LEMMA A.5. *Under* (C0) *and* (C2$'$) *for* $G = F$ *and* $F_n$ *and very large* $n$:

(a) $\sup_{x \in \mathbb{R}^d} |PD(x, F_n) - PD(x, F)| = o(1)$, *a.s.*,
(b) $PD(x, G)$ *is Lipschitz continuous in* $x \in \mathbb{R}^d$, *a.s.*,
(c) $\partial PD^\alpha(G) = \{x : PD(x, G) = \alpha\}$, *a.s.*,
(d) $PD^{(\alpha+\eta)}(F) \subseteq PD^{(\alpha+\eta/2)}(F_n) \subseteq PD^\alpha(F)$ *a.s. for any* $0 < \eta < \alpha^* - \alpha$.

Now we prove the theorem. Condition (C1$'$) implies (C0). By Lemma A.5, $PD^\alpha(F_n)$ is nonempty and contains the origin a.s. for large $n$. Hence $R^\alpha(u, F_n)$ is well defined a.s. Condition (C1$'$) also implies that there is a unit vector $v(x)$ such that $g(x, v(x), F) = O(x, F)$ for any $x \in \mathbb{R}^d$. Let $v(u) := v(R(u)u)$. Likewise, we have a unit vector $v_n(u) := v(R_n(u)u)$. By virtue of Lemma A.5, $O(R(u)u, F) = O(R_n(u)u, F_n) = \beta(\alpha)$ a.s. for sufficiently large $n$. Hence for any $v \in S^{d-1}$

$$R(u)u'v(u) = \beta(\alpha)\sigma(F_{v(u)}) + \mu(F_{v(u)}), \qquad R(u)u'v \leq \beta(\alpha)\sigma(F_v) + \mu(F_v).$$

Likewise, we can have the same displays for $R_n(u)$, $v_n(u)$ and $v$. These give

$$(\beta(\alpha)[\sigma(F_{nv_n(u)}) - \sigma(F_{v_n(u)})] + (\mu(F_{nv_n(u)}) - \mu(F_{v_n(u)})))/u'v_n(u)$$
$$\leq R_n(u) - R(u)$$
$$\leq (\beta(\alpha)[\sigma(F_{nv(u)}) - \sigma(F_{v(u)})] + (\mu(F_{nv(u)}) - \mu(F_{v(u)})))/u'v(u).$$

If we can show that $\inf_{u \in S^{d-1}} |u'v(u)| > 0$ and $\inf_{u \in S^{d-1}} |u'v_n(u)| > 0$ almost surely for large $n$, then the theorem follows in a straightforward fashion from (C2′).

The proof for $\inf_{u \in S^{d-1}} |u'v(u)| > 0$ is given in the proof of Lemma A.2. The argument for proving $\inf_{u \in S^{d-1}} |u'v_n(u)| > 0$ a.s. for sufficiently large $n$ is the same. But we need the following two almost sure results for sufficiently large $n$:

(C0″) $\sup_{\|u\|=1} \mu(F_{nu}) < \infty$, $0 < \inf_{\|u\|=1} \sigma(F_{nu}) \leq \sup_{\|u\|=1} \sigma(F_{nu}) < \infty$,

and $O(0, F_n) < O(R_n(u)u, F_n)$ a.s. The first one (C0″) follows from (C0) and (C2′). The second one follows from Lemma A.5 since the origin is an interior point of $PD^{\alpha+\delta}(F) \subset PD^{\alpha}(F_n)$ a.s. for some $0 < \delta < \alpha^* - \alpha$ and sufficiently large $n$. □

PROOF OF THEOREM 5. The following lemma about the continuity of $R^\alpha(u, F)$ in $u \in S^{d-1}$ is needed in the proof of the theorem.

LEMMA A.6. *Under* (C0) *and* (C2′), $R^\alpha(u, F_n)$ *is continuous in $u$ for large $n$.*

We now prove Theorem 5. Following the proof of Theorem 4, we have

$$(\beta(\alpha)[\sigma(F_{nv_n(u)}) - \sigma(F_{v_n(u)})] + (\mu(F_{nv_n(u)}) - \mu(F_{v_n(u)})))/u'v_n(u)$$
$$\leq R_n(u) - R(u)$$
$$\leq (\beta(\alpha)[\sigma(F_{nv(u)}) - \sigma(F_{v(u)})] + (\mu(F_{nv(u)}) - \mu(F_{v(u)})))/u'v(u),$$

and $\inf_{u \in S^{d-1}} u'v(u) > 0$ and $\inf_{u \in S^{d-1}} u'v_n(u) > 0$ almost surely for $n$ large.

By the compactness of $S^{d-1}$, the continuity in (C1′) is uniform in $u \in S^{d-1}$. This, in conjunction with the last display, (C3′) and standard results on empirical processes (see, e.g., Problem II.18, Approximation Lemma II.25, Lemma II.36, Equicontinuity Lemma VII.15, and (the central limit theorem for empirical processes) Theorem VII.21 of [20], or see [32]), gives the desired results if we can show that $v_n(u) \to v(u)$ uniformly in the unit vector $u$ as $n \to \infty$. The latter can be done in much the same way as the uniform convergence of $v_\varepsilon(u) \to v(u)$ as $\varepsilon \to 0$ in the proof of Theorem 1. □

PROOF OF THEOREM 6. The desired result follows if we show that the numerator and the denominator of $PTM(F_n)$ converge a.s. to those of $PTM(F)$, respectively. Clearly it suffices to treat just the numerator. The given conditions imply (C0), which, combining with (C2'), Lemma A.5 and the continuity of $w^{(1)}$, yields

$$w(PD(x, F_n)) = w(PD(x, F)) + o(1) \qquad \text{a.s. and uniformly in } x \in \mathbb{R}^d. \tag{23}$$

The boundedness of $PD^{\alpha'}(F)$ for any $\alpha' > 0$ (see Theorem 2.3 of [33]) and Lemma A.5 imply the almost sure boundedness of $PD^{\alpha}(F_n)$ for sufficiently large $n$. This, together with (23), implies that

$$\int_{PD(x,F_n) \geq \alpha} w(PD(x, F_n)) x \, dF_n(x)$$
$$= \int_{PD(x,F_n) \geq \alpha} w(PD(x, F)) x \, dF_n(x) + o(1) \qquad \text{a.s.}$$

The desired result follows if we can show that

$$\int_{PD(x,F_n) \geq \alpha} w(PD(x, F)) x \, dF_n(x)$$
$$- \int_{PD(x,F) \geq \alpha} w(PD(x, F)) x \, dF(x) = o(1) \qquad \text{a.s.}$$

In light of the (a.s.) compactness of $PD^{\alpha}(F)$ and $PD^{\alpha}(F_n)$, Lemma A.5 and Lebesgue's dominated convergence theorem, we see that

$$\int [I(PD(x, F_n) \geq \alpha) - I(PD(x, F) \geq \alpha)] w(PD(x, F)) x \, dF(x) = o(1) \qquad \text{a.s.}$$

Thus we only need to show that

$$\int I(PD(x, F_n) \geq \alpha) w(PD(x, F)) x \, d(F_n - F)(x) = o(1) \qquad \text{a.s.} \tag{24}$$

Let $\delta \in (0, \alpha)$. Then $PD^{\alpha-\delta}(F)$ is convex and compact (Theorem 2.3 of [33]). By Lemma A.5, $PD^{\alpha}(F_n) \subset PD^{\alpha-\delta}(F)$ a.s. for sufficiently large $n$. This and the convexity of $O(\cdot, F_n)$ imply that $PD^{\alpha}(F_n)$ is convex and compact and contained in $PD^{\alpha-\delta}(F)$ a.s. Define $\mathcal{C} = \{C : C \subset PD^{\alpha-\delta}(F) \text{ is compact and convex}\}$. Then $PD^{\alpha}(F_n) \in \mathcal{C}$ a.s. for sufficiently large $n$. By a well-known result of Ranga Rao ([22], Theorem 4.2), $\mathcal{C}$ is an F-Glivenko–Cantelli class (see [32], for the corresponding definition and related discussion) and (24) follows from the boundedness of $w$ and $PD^{\alpha}(F_n)$ (a.s.). □

PROOF OF THEOREM 7. The following representation of $H_n(x) := \sqrt{n} \times (PD(x, F_n) - PD(x, F))$, established in Lemma 5.2 of [35], is needed.

LEMMA A.7. *Let $\mu(F_u)$ and $\sigma(F_u)$ be continuous in $u$ and $\sigma(F_u) > 0$ for $u \in S^{d-1}$. Then under (C3′) and (C4) we have for $\nu_n = \sqrt{n}(F_n - F)$*

$$\sqrt{n}(PD(x, F_n) - PD(x, F)) = \int \tilde{h}(y, x)\nu_n(dy) + o_p(1) \quad \text{uniformly in } x \in B$$

*with $\tilde{h}(y, x) = (O(x, F)f_2(y, u(x)) + f_1(y, u(x)))/(\sigma(F_{u(x)})(1 + O(x, F))^2)$.*

We now prove Theorem 7. First we note that

$$PTM(F_n) - PTM(F) = \frac{\int_{PD^\alpha(F_n)} (x - PTM(F))w(PD(x, F_n))\, dF_n(x)}{\int_{PD^\alpha(F_n)} w(PD(x, F_n))\, dF_n(x)}.$$

Following the proof of Theorem 6, we can see immediately that

$$(25) \quad \int_{PD^\alpha(F_n)} w(PD(x, F_n))\, dF_n(x) = \int_{PD^\alpha(F)} w(PD(x, F))\, dF(x) + o(1) \quad \text{a.s.}$$

Decompose the numerator of $PTM(F_n) - PTM(F)$ into three parts,

$$I_{1n} = \int_{PD^\alpha(F_n)} x^* W(PD(x, F_n))\, dF_n(x) - \int_{PD^\alpha(F)} x^* W(PD(x, F_n))\, dF_n(x),$$

$$I_{2n} = \int_{PD^\alpha(F)} x^* W(PD(x, F_n))\, dF_n(x) - \int_{PD^\alpha(F)} x^* W(PD(x, F))\, dF_n(x),$$

$$I_{3n} = \int_{PD^\alpha(F)} x^* W(PD(x, F))\, d(F_n(x) - F(x)),$$

where $x^* = x - PTM(F)$. Obviously

$$(26) \quad \sqrt{n} I_{3n} = \int (x - PTM(F))W(PD(x, F))I(x \in PD^\alpha(F))\, d\nu_n(x).$$

We now work on $I_{2n}$. By (C3′) and the central limit theorem for empirical processes (see, e.g., Theorem VII. 21 of [20]) we have that

$$(27) \quad \begin{aligned} \sup_{\|u\|=1} |\mu(F_{nu}) - \mu(F_u)| &= O_p(1/\sqrt{n}), \\ \sup_{\|u\|=1} |\sigma(F_{nu}) - \sigma(F_u)| &= O_p(1/\sqrt{n}), \end{aligned}$$

which then imply (see Theorem 2.2 and Remark 2.5 of [33]) that

$$\sup_{x \in \mathbb{R}^d} (1 + \|x\|)|PD(x, F_n) - PD(x, F)| = O_p(1/\sqrt{n}). \tag{28}$$

By the continuous differentiability of $w$ on $[0, 1]$ and Lemma A.5, we have that

$$w(PD(x, F_n)) - w(PD(x, F)) \\ = (w^{(1)}(PD(x, F)) + o(1))H_n(x)/\sqrt{n} \quad \text{a.s.} \tag{29}$$

uniformly in $x \in \mathbb{R}^d$. This, the boundedness of $PD^\alpha(F)$ and (28) imply

$$\sqrt{n} I_{2n} = \int_{PD^\alpha(F)} (x - PTM(F)) w^{(1)}(PD(x, F)) H_n(x) \, dF_n(x) + o_p(1). \tag{30}$$

We now show that

$$\sqrt{n} I_{2n} = \int_{PD^\alpha(F)} (x - PTM(F)) w^{(1)}(PD(x, F)) H_n(x) \, dF(x) + o_p(1). \tag{31}$$

Clearly, we can view $H_n(\cdot)$ for every $n$ as a map into $l^\infty(\mathbb{R}^d)$, the space of all uniformly bounded, real functions on $\mathbb{R}^d$. By Lemma A.7 (and its proof; see [35]) and Theorem 1.5.4 of [32], we see that $H_n$ is asymptotically tight on $B$ (the set defined in the theorem). Consequently, for every $\varepsilon > 0$ there are finitely many continuous functions $h_1, \ldots, h_k$ on $B$ such that

$$\limsup_{n \to \infty} P \left\{ \min_{1 \leq i \leq k} \|H_n - h_i\|_\infty > \varepsilon \right\} \leq \varepsilon.$$

Since the functions $I(PD(x, F) \geq \alpha)(x - PTM(F))w^{(1)}(PD(x, F))h_i(x)$ are bounded and continuous for $x \in \{PD(x, F) \geq \alpha\}$, hence

$$\left\| \int_{PD^\alpha(F)} (x - PTM(F)) w^{(1)}(PD(x, F)) H_n(x) \, d(F_n - F)(x) \right\|$$

$$\leq \max_{1 \leq i \leq k} \left\| \int_{PD^\alpha(F)} (x - PTM(F)) w^{(1)}(PD(x, F)) h_i(x) \, d(F_n - F)(x) \right\|$$

$$+ (2\varepsilon) \sup_{x \in PD^\alpha(F)} \|(x - PTM(F)) w^{(1)}(PD(x, F))\|$$

$$= O(\varepsilon) + o(1)$$

with asymptotic probability not less than $1 - \varepsilon$. Thus we obtain (31), which, in conjunction with Lemma A.7, (C3′) and Fubini's theorem, gives

$$\sqrt{n} I_{2n} = \int \left( \int_{PD^\alpha(F)} (y - PTM(F)) w^{(1)}(PD(y, F)) h'(x, y) \, dF(y) \right) \\ \times d\nu_n(x) + o_p(1). \tag{32}$$

We now work on $I_{1n}$. Let $\Delta_n(x) = I(PD(x, F_n) \geq \alpha) - I(PD(x, F) \geq \alpha)$. By (28) and (29) and the boundedness of $PD^\alpha(F)$ and $PD^\alpha(F_n)$ (a.s. for sufficiently large $n$) (see Lemma A.5), we have

$$\sqrt{n} I_{1n} = \sqrt{n} \int \Delta_n(x)(x - PTM(F)) w(PD(x, F)) \, dF_n(x)$$
$$+ \int \Delta_n(x)(x - PTM(F)) w^{(1)}(PD(x, F)) H_n(x) \, dF_n(x) + o_p(1),$$

for sufficiently large $n$. Call two terms on the right-hand side $I_{1n1}$ and $I_{1n2}$, respectively.

We first show that $I_{1n2} = o_p(1)$. Observe that by (28) we have

$$\|I_{1n2}\| \leq |O_p(1)| \int \|x - PTM(F)\| |\Delta_n(x) w^{(1)}(PD(x, F))| \, dF_n(x).$$

Invoking the Skorohod (representation) theorem, we assume that $Y_n$ and $Y$ are defined on the probability space $(\Omega, \mathcal{F}, P)$ such that $Y_n - Y = o(1)$ a.s. $(P)$ and $F_{Y_n} = F_n$ and $F_Y = F$. By changing the variables in the above integral, we have

$$\int \|x^*\| |\Delta_n(x) w^{(1)}(PD(x, F))| \, dF_n(x)$$
$$= \int_\Omega \|Y_n^*\| |\Delta_n(Y_n) w^{(1)}(PD(Y_n, F))| \, dP,$$

where $Y_n^* = Y_n - PTM(F)$ and $\Delta_n(Y_n) \to 0$ a.s. by Lemma A.5. This, Lemma A.5 and Lebesgue's dominated convergence theorem yield immediately $I_{1n2} = o_p(1)$.

We now show that

$$(33) \quad I_{1n1} = \sqrt{n} \int \Delta_n(x)(x - PTM(F)) w(PD(x, F)) \, dF(x) + o_p(1).$$

This can be accomplished by utilizing the results of an F-Donsker class of functions and the asymptotic equicontinuity (see [32]) and the fact $\int (I(PD(x, F_n) \geq \alpha) - I(PD(x, F) \geq \alpha))^2 \, dF(x) \to 0$. Observe that

$$\sqrt{n} \int \Delta_n(x)(x - PTM(F)) w(PD(x, F)) \, dF(x)$$
$$= \sqrt{n} \int_{S^{d-1}} \left[ \int_{R^\alpha(u,F)}^{R^\alpha(u,F_n)} (ru - PTM(F)) w(PD(ru, F)) f(ru) |J(u, r)| \, dr \right] du.$$

Let $\theta_n(u)$ be a point in between $R^\alpha(u, F_n)$ and $\mathbb{R}^\alpha(u, F)$. Then by Theorem 4, $J(u, \theta_n(u)) = J(u, R^\alpha(u, F)) + o(1)$ a.s. uniformly in $u \in S^{d-1}$. By Theorem 4 and Lemma A.5, $w(PD(\theta_n(u)u, F)) = w(PD(R^\alpha(u, F)u, F)) + o(1) = w(\alpha) + o(1)$ a.s. and uniformly in $u \in S^{d-1}$. Finally by the continuity of $f$ [in a small neighborhood of $\partial PD^\alpha(F)$] and of $R^\alpha(u, F_n)$ and $R^\alpha(u, F)$ uniformly in $u$ (see Lemmas A.3 and A.6) for large $n$, the compactness of $S^{d-1}$ and Theorem 4, $f(\theta_n(u)u) = f(R^\alpha(u,F)u) + o(1)$, a.s. uniformly in $u \in S^{d-1}$ for large $n$. These, (33), the preceding display, the mean value theorem, the uniform continuity in $u$ of $w(PD(R^\alpha(u,F)u,F))$ and $J(u,R^\alpha(u,F))$, Theorems 4 and 5, and Fubini's theorem, yield

$$I_{1n1} = \int_{S^{d-1}} (R^*(u)w(\alpha)|J(u,R(u))|f(R(u)u) + o(1))K_n(u)\,du + o_p(1)$$

$$= \int \left[\int_{S^{d-1}} R^*(u)w(\alpha)|J(u,R(u))|f(R(F)u)k(x,R(u)u)\,du\right]d\nu_n(x) + o_p(1),$$

for $R^*(u) = R(u)u - PTM(F)$, $K_n(u) := \sqrt{n}(R(u,F_n) - R(u,F))$. This gives

$$\sqrt{n}I_{1n} = \int \left[\int_{S^{d-1}} R^*(u)w(\alpha)|J(u,R(u))|f(R(u)u)k(x,R(u)u)\,du\right]d\nu_n(x) + o_p(1),$$

which, combining with (32), (26) and (25), gives the desired result. □

**Acknowledgments.** The author thanks the referees, the Associate Editor and the Editor Morris L. Eaton for insightful comments and suggestions.

DEPARTMENT OF STATISTICS AND PROBABILITY
MICHIGAN STATE UNIVERSITY
EAST LANSING, MICHIGAN 48824
USA
E-MAIL: zuo@stt.msu.edu